\RequirePackage{fix-cm}
\documentclass[smallextended]{svjour3}       
\smartqed  
\usepackage{graphicx}
\usepackage{color}
\RequirePackage{amsmath}%
\usepackage{url}%
\usepackage{amsfonts}%

\renewcommand{\Re}{\mathop{\mathrm{Re}}}
\renewcommand{\Im}{\mathop{\mathrm{Im}}}
\renewcommand{\i}{\mathrm{i}}
\newcommand{\bI}{{\bf I}}

\newcommand{\bM}{{\bf M}}
\newcommand{\bN}{{\bf N}}

\newcommand{\C}{{\mathbb C}}
\newcommand{\R}{{\mathbb R}}
\newcommand{\D}{{\mathbb D}}

\newcommand{\capa}{\mathrm{cap}}
\newcommand{\K}{\mathcal{K}}

\journalname{ }

\begin{document}

\title{Numerical computation of a preimage domain for an infinite strip with rectilinear slits}

\titlerunning{Numerical computation of a preimage domain}        

\author{El Mostafa Kalmoun          \and
        Mohamed M. S. Nasser        \and 
        Matti Vuorinen }

\authorrunning{E. Kalmoun \and M. Nasser \and M. Vuorinen } 

\institute{E. Kalmoun \at
              School of Science and Engineering, Al Akhawayn University in Ifrane, PO Box 104, Ifrane 53000, Morocco \\
              \email{E.Kalmoun@aui.ma} 
		 \and
		M. Nasser \at
              Department of Mathematics, Statistics \& Physics, Wichita State University, Wichita, KS 67260-0033, USA \\
              \email{mms.nasser@wichita.edu}           
        \and
        M. Vuorinen \at
              Department of Mathematics and Statistics, University of Turku, FI-20014 Turku, Finland \\
              \email{vuorinen@utu.fi} 
}

\date{Received: date / Accepted: date}

\maketitle

\begin{abstract}
Let $\Omega$ be the multiply connected domain in the extended complex plane $\overline{{\mathbb C}}$ obtained by removing $m$ non-overlapping rectilinear segments from the infinite strip $S=\{z\,:\, \left|\Im z\right|<\pi/2\}$. 
In this paper, we present an iterative method for numerical computation of a conformally equivalent bounded multiply connected domain $G$ in the interior of the unit disk $\D$ and the exterior of $m$ non-overlapping smooth Jordan curves. We demonstrate the utility of the proposed method through two applications. First, we estimate the capacity of condensers of the form $(S,E)$ where $E\subset S$ be a union of disjoint segments. Second, we determine the streamlines associated with uniform incompressible, inviscid and irrotational flow past disjoint segments in the strip $S$.
\end{abstract}

\section{Introduction}\label{section1}

Let $\Omega$ be the multiply connected domain in the extended complex plane $\overline{\C}=\C\cup\{\infty\}$ obtained by removing $m$ non-overlapping rectilinear segments from the infinite strip $S=\{z\,:\, \left|\Im z\right|<\pi/2\}$. 
Solving boundary value problems in such a domain with complicated boundaries is not as easy as it is for domains with smooth boundaries. 
A possible remedy is to find a conformal mapping from $\Omega$ onto a multiply connected domain $G$ bordered by smooth Jordan curves. 

The above domain $\Omega$ is one of the canonical domains for conformal mapping of multiply connected domains~\cite[p.~128]{wen}.
An efficient method for numerical computation of the conformal mapping $\Phi$ from domains with smooth boundaries $G$ onto the canonical domain $\Omega$ is presented in~\cite{NF13}. Still, in this method, the domain $G$ is supposed given while $\Omega$ is unknown and should be computed alongside the conformal mapping $\Phi$ from $G$ onto $\Omega$. 
On the contrary, in this paper, we take up the case when $\Omega$ is known. Our objective will be then to find an unknown preimage domain $G$ bordered by smooth Jordan curves, and to determine a conformal mapping $\Phi$ from $G$ onto $\Omega$. This means that the method presented in~\cite{NF13} is not directly applicable in our context, but this method is still useful to develop an iterative scheme for finding the unknown preimage domain $G$ as well as the conformal mapping $\Phi$.  
In this way, the inverse mapping $\Phi^{-1}$ would be the desired conformal mapping from the domain $\Omega$ onto the domain $G$.
It is known that Laplace equation in the plane is invariant under conformal mappings. Thus, with the help of the inverse conformal mapping $\Phi^{-1}$, boundary value problems for the Laplace equation in the domain $\Omega$ will be solved easily in the domain $G$.

The idea of our iterative method is similar in spirit to the approach employed in~\cite{AST13,N19,NG18} but for other canonical domains. The proposed numerical method will be useful in solving several problems in strip with rectilinear slit domains. Applications of the proposed method to two types of such problems will be considered in this paper. 

As a main illustration of how useful is the developed numerical conformal mapping method of this paper, we study a domain functional, the conformal capacity of a condenser.
Recall that a \emph{condenser} is a pair $(D,E)$ where $ D \subset
\C$ is a domain and $E \subset D$ is a non-empty compact set, and the  
\emph{conformal capacity} is defined as follows
\begin{align}\label{def_condensercap}
{\rm cap}(D,E)=\inf_{u\in \mathcal{A}}\int_{D}|\nabla u|^2 \, dm,
\end{align}
where $dm$ stands for the 2-dimensional Lebesgue measure and $\mathcal{A}$ is the family of all harmonic functions in $D$ with values in $[1,\infty)$ on $E$ and approaching $0$ on $\partial D$.  
If $D\setminus E$ is a multiply connected domain such that all of its boundary components are smooth Jordan curves, then the infimum in \eqref{def_condensercap} is known to be attained by a harmonic function \cite[p. 65]{ah}. In fact, 
this extremal function is the solution of the 
Laplace equation with boundary values equal to $1$ on $E$ and vanishing on $\partial D$. The conformal capacity, which is a mathematical model
of the capacity of a physical condenser, has numerous applications also 
to potential theory and to geometric function theory \cite{du,hkv,sol}.

As is generally known, the set function $E \mapsto {\rm cap}(D,E)$
looks similar to an outer measure \cite[Lemma 7.1, Theorem 9.6]{hkv}. 
For instance, the subadditivity property of the condenser capacity
for compact sets $E_j \subset D, j=1,...,m,$ says that
\begin{equation} \label{subadd}
{\rm cap} (D, \cup_{j=1}^m E_j)\le \sum_{j=1}^m {\rm cap} (D,E_j).
\end{equation}
This is reminiscent of the subadditivity
of  Lebesgue measurable sets $E_j \subset \R,
j=1,2,...,$
\[
{\rm m}( \cup_{j=1}^{\infty}{E_j}) \le \sum_{j=1}^{\infty}{\rm m}({E_j}) 
\]
with equality for separate sets. For the capacity case \eqref{subadd} the
situation is different, no equality statement is known for separate sets. Moreover, even in the
case $m=2$ it is not easy to give non-trivial examples of sets with strict
inequality in \eqref{subadd}.

Our numerical computation leads
to two novel observations. First, we 
will show here that an asymptotic equality holds,
the lower bound in \eqref{subadd} will be arbitrarily close to the upper
bound for some sets, far away from each other. This asymptotic equality is
a  manifestation of some kind of  ``\emph{weak additivity}'' for such sets. Second,
our computational experiments have revealed an inequality for the capacity of condensers in the strip domain $S.$ More precisely, if the condenser is of the form
$(S,E)$ and $E= E_1 \cup E_2 \subset {\mathbb R}$, $E_1$ and $E_2$ are segments, then
\begin{equation}\label{eq:ineq}
 {\rm cap}(S,E)\ge {\rm cap}(S,H)
\end{equation}
where $H \subset {\mathbb R}$ is a segment with diameter equal to the sum 
of diameters of $E_1$ and $E_2.$ Due to the conformal invariance of the capacity,
the inequality \eqref{eq:ineq} admits various extensions to simply connected plane
domains, for instance to the case when the strip domain $S$ is replaced by the
unit disk. In that case we have to use the hyperbolic metric. The case of multiply
connected domains seems to offer problems for further research.

In the second application, we consider computing the complex potential for a uniform inviscid and incompressible flow past multiple disjoint segment obstacles in the strip $S$ in the case when the circulations around the segments are zeros. 
A numerical method for approximating such complex potentials when the obstacles have smooth Jordan curves boundaries has been presented in~\cite{sak12} where the domain $\Omega$ was called a channel domain.
The idea used in~\cite{sak12} is based on constructing a conformal mapping $w=F(z)$ from $\Omega$ onto the domain obtained by removing horizontal slits from the strip $S$. Then, $W(z)=F(z)$ represents a complex potential for the uniform flow in $\Omega$. The streamlines for a uniform flow in $\Omega$ are then the contour plots of the imaginary part of the complex potential $W(z)$.

In this paper, the geometry of the domain $\Omega$ is more complicated as the obstacles are slits. However, the method used here is similar to the method presented in~\cite{sak12} and hinges on the numerical computation of the conformal mapping $w=F(z)$ from $\Omega$ onto a domain $H$ obtained by removing horizontal slits from the strip $S$. To compute such a conformal mapping, we first apply the proposed iterative method to compute a conformal mapping from $\Omega$ to a domain $G$ bordered by smooth Jordan curves, and then use the method of~\cite{NF13} to conformally map $G$ onto the domain $H$ obtained by removing horizontal slits from the strip $S$.

Further applications of the proposed iterative method are possible. For example, with the help of the iterative method, one can extend the method presented in~\cite{NasVla} for simulating local fields in carbon nanotube (CNT) reinforced composites for infinite strip with circular voids and elliptic CNTs to the case of slit CNTs. Another example is to extend the method in~\cite{NG18} to compute ideal fluid flow in channel domains with multiple slit stirrers. 

It is worth mentioning that the above multiply connected domain $\Omega$ is a degenerate case of a multiply connected polygonal domain. 
For simply connected domains, the Schwarz–Christoffel formula provides us with an explicit formula for computing the conformal mapping from the unit disk onto a given polygonal domain. 
Recently, the Schwarz--Christoffel formula has been generalized to compute conformal mappings from multiply connected circular domains onto multiply connected polygonal domains by DeLillo, Elcrat, and Pfaltzgraff~\cite{dep} and Crowdy~\cite{crow-05,crow-07} (see also~\cite{crow-bad,crow-20}). 
However, to use such explicit formulas, we need to solve systems of non-linear equations to determine the preimages of the vertices of the polygons and, for multiply connected domains, the centers and the radii of the circles.

\section{Conformally mapping  a strip with rectilinear slits onto a domain with smooth Jordan boundaries}
\label{sec:strip}

Suppose that $\Omega$ is the canonical multiply connected domain obtained by removing $m$ non-overlapping rectilinear slits $[a_j,b_j]$ from the infinite strip $S=\{z\,:\, \left|\Im z\right|<\pi/2\}$ where $a_j,b_j\in S$ are complex numbers, $j=1,\ldots,m$, i.e. $\Omega=S \setminus \cup_{j=1}^{m}[a_j,b_j]$; see Figure~\ref{fig1} for an example of $\Omega $ with $m=4$.

\subsection{Boundary integral equation for the conformal mapping}

In this subsection, we briefly review the method presented in~\cite{NF13} to compute a conformal mapping from  a multiply connected domain $G$ with smooth Jordan curves as boundaries onto the above domain $\Omega$. The domain $G$ is of the form $\mathbb{D}\setminus \cup_{j=1}^m \overline{E_j}$ where the $\{E_j : j-1,\ldots m\}$ is a collection of disjoint simply connected  domains that are bounded by smooth Jordan curves $\Gamma_1,...,\Gamma_m$. 
The external boundary of $G$ is the unit circle $\Gamma_0 := \partial \mathbb{D}.$
In this setting, we have the existence of a unique conformal mapping $\Phi$ from $G$ onto $\Omega$ such that~\cite[p.~128]{wen}
\begin{equation}\label{eq:strip-cond}
\Phi(\pm 1)=\pm\infty+\i 0 \quad\mbox{and}\quad \Phi(\i)=\frac{\pi}{2}\i.
\end{equation}
For $j=1,\ldots,m$, let $\ell_j=|b_j-a_j|$ be the length of the segment $L_j=[a_j,b_j]$,  $c_j=(a_j+b_j)/2$ be the center of $L_j$, and  $\theta_j$ be the angle between $L_j$ and the positive real axis. The values of the real constants $\ell_1,\ldots,\ell_m$ and the complex constants $c_1,\ldots,c_m$ are undetermined and should be computed alongside the conformal mapping $\Phi$. These constants are uniquely determined by the domain $G$. On the other hand, the values of the angles $\theta_1,\ldots,\theta_m$ can be fixed in advance. 

We parametrize each boundary component $\Gamma_j$ by a $2\pi$-periodic complex-valued function $\eta_j(t)$, $t\in J_j:=[0,2\pi]$, $j=0,\ldots,m$. 
Henceforth, we define the total parameter domain $J$ as the disjoint union of the $m+1$ intervals $J_j=[0,2\pi]$, $j=0,\ldots,m$. In this way, 
the whole boundary $\Gamma$ is parametrized by
\begin{equation}\label{eq:eta}
\eta(t)= \left\{ \begin{array}{l@{\hspace{0.5cm}}l}
\eta_0(t),&t\in J_0,\\
\vdots & \\
\eta_m(t),&t\in J_m.
\end{array}
\right.
\end{equation}
See~\cite{Nas-ETNA} for more details.

We take $\theta_0=0$, assume that $\theta_1,\ldots,\theta_m$ are given real constants, and consider the complex function $A$ defined by
\begin{equation}\label{eq:A}
A(t) = e^{\i(\frac{\pi}{2}-\theta(t))}(\eta(t)-\alpha),
\end{equation}
where $\alpha\in G$ is given and $\theta$ is defined on $J$ by 
\[
\theta(t)=\theta_j \quad {\rm for}\quad t\in J_j, \quad j=0,\ldots,m,
\]
i.e., the function $\theta$ is constant on each interval $J_j$.
Then, the kernel $N(s,t)$ defined on $J\times J$ by 
\begin{equation}\label{eq:N}
N(s,t) :=
\frac{1}{\pi}\Im\left(\frac{A(s)}{A(t)}\frac{\dot\eta(t)}{\eta(t)-\eta(s)}\right) \quad {\rm for}\quad (s,t)\in J\times J,
\end{equation}
is known as the \emph{generalized Neumann kernel}.
We also define the kernel $M(s,t)$ on $J\times J$ by 
\begin{equation}\label{eq:M}
M(s,t) :=
\frac{1}{\pi}\Re\left(\frac{A(s)}{A(t)}\frac{\dot\eta(t)}{\eta(t)-\eta(s)}\right)  \quad {\rm for}\quad (s,t)\in J\times J.
\end{equation}
Note that $N(s,t)$ is continuous while $M(s,t)$ is singular with its singular part involving the cotangent function. Hence, the integral operator $\bN$ with the kernel $N(s,t)$ is compact and the integral operator $\bM$ with the kernel $M(s,t)$ is singular. Further details can be found in~\cite{Weg-Nas}.

The method proposed in~\cite{NF13} for computing the conformal mapping $\Phi$ from the domain $G$ onto the domain $\Omega$ is summarized in the following theorem. For this, let us set
\begin{equation*}\label{eq:gama-strip}
	\gamma(t)= \left\{ \begin{array}{l@{\hspace{0.5cm}}l}
		0,&t\in J_0,\\
		\Im\left[e^{-\i\theta_j}\Psi(\eta(t))\right],&t\in J_j,\quad j=1,\ldots,m, \\
	\end{array}
	\right.
\end{equation*}
where 
\begin{equation}\label{eq:Psi}
	\Psi(w) = \log\frac{1+w}{1-w}.
\end{equation}
The function $\zeta=\Psi(w)$ is a conformal mapping from the unit disk $\D$ onto the infinite strip $-\pi/2<\Im\zeta<\pi/2$.

\begin{theorem}[\cite{NF13}]\label{thm:cm-strip}
If $\rho$ is the unique solution of the boundary integral equation
\begin{equation}\label{eq:ie}
(\bI-\bN)\rho=-\bM\gamma,
\end{equation}
and the piecewise constant function $h$ is given by
\begin{equation}\label{eq:h}
h=[\bM\rho-(\bI-\bN)\gamma]/2,
\end{equation}
then the conformal mapping $\Phi$ from $G$ onto $\Omega$ is given by
\begin{equation}
\label{eq:Phi-f}
\Phi(w)=-\i f(\i)+w f(w)+\Psi(w), \quad w\in G\cup\Gamma,
\end{equation}
where $f$ is the analytic function in $G$ with the boundary values
\begin{equation}\label{eq:strip-f}
A(t)f(\eta(t))=\gamma(t)+h(t)+\i\rho(t).
\end{equation}
\end{theorem}

\subsection{Computing the preimage domain $G$}\label{sec:itr}

We assume that the strip with rectilinear slits domain $\Omega$ is given. This means that the values of the constants $\ell_j$, $c_j$, and $\theta_j$ are now known for $j=1,\ldots,m$.
The method shown in Theorem~\ref{thm:cm-strip} will be used in this subsection to develop an iterative scheme to find a bounded multiply connected preimage domain $G$ in the exterior of $m$ smooth Jordan curves $\Gamma_1,\ldots,\Gamma_m$ and the interior of the unit circle $\Gamma_0$ as well as the conformal map $z=\Phi(w)$ from $G$ onto $\Omega$ that satisfies the normalization conditions~\eqref{eq:strip-cond}.

\begin{figure}[ht]\label{fig1}
\centerline{
\scalebox{0.3}{\includegraphics[trim=0 -1.0cm 0 0,clip]{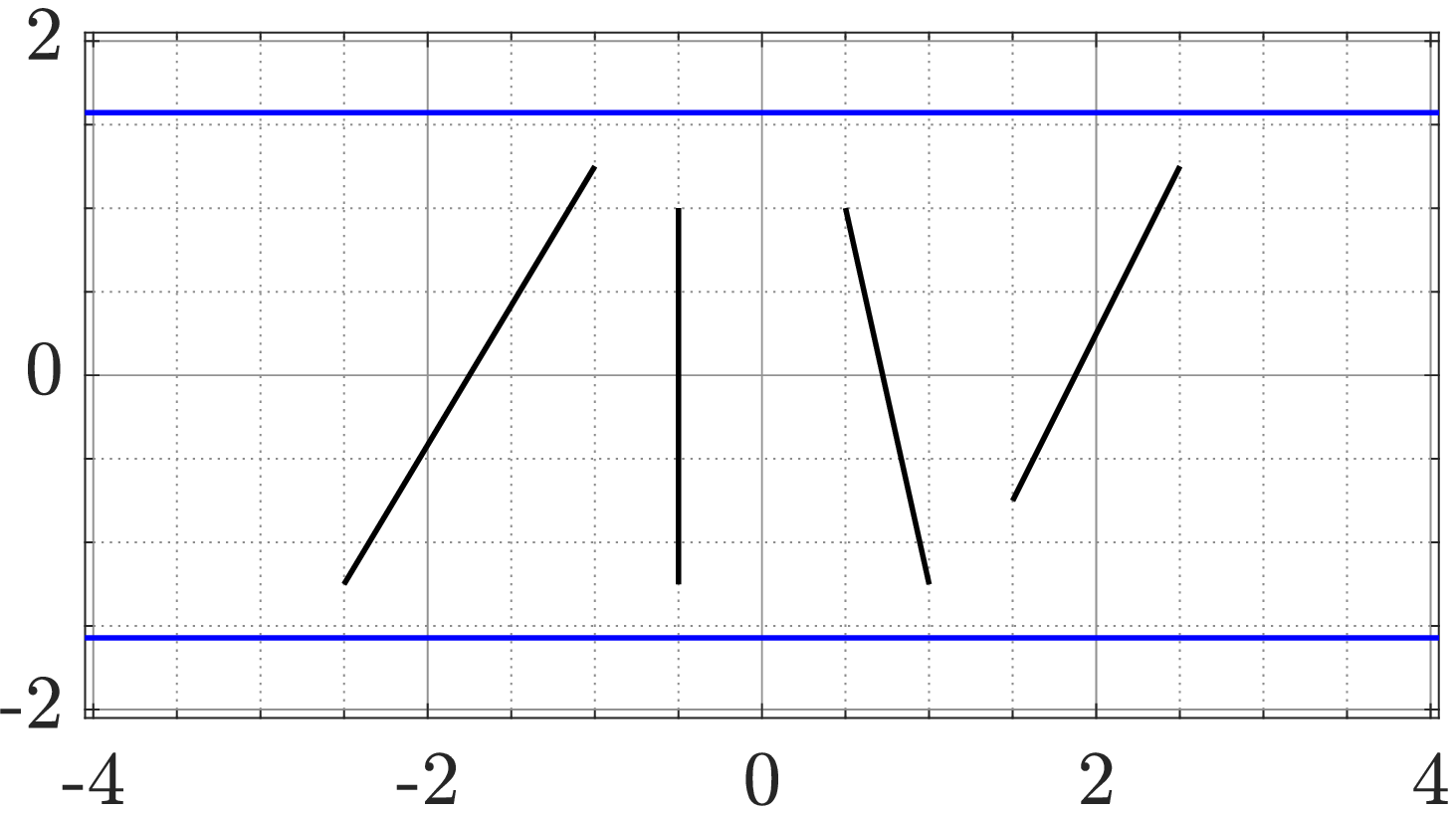}}
\hfill
\scalebox{0.3}{\includegraphics[trim=0 -1.0cm 0 0,clip]{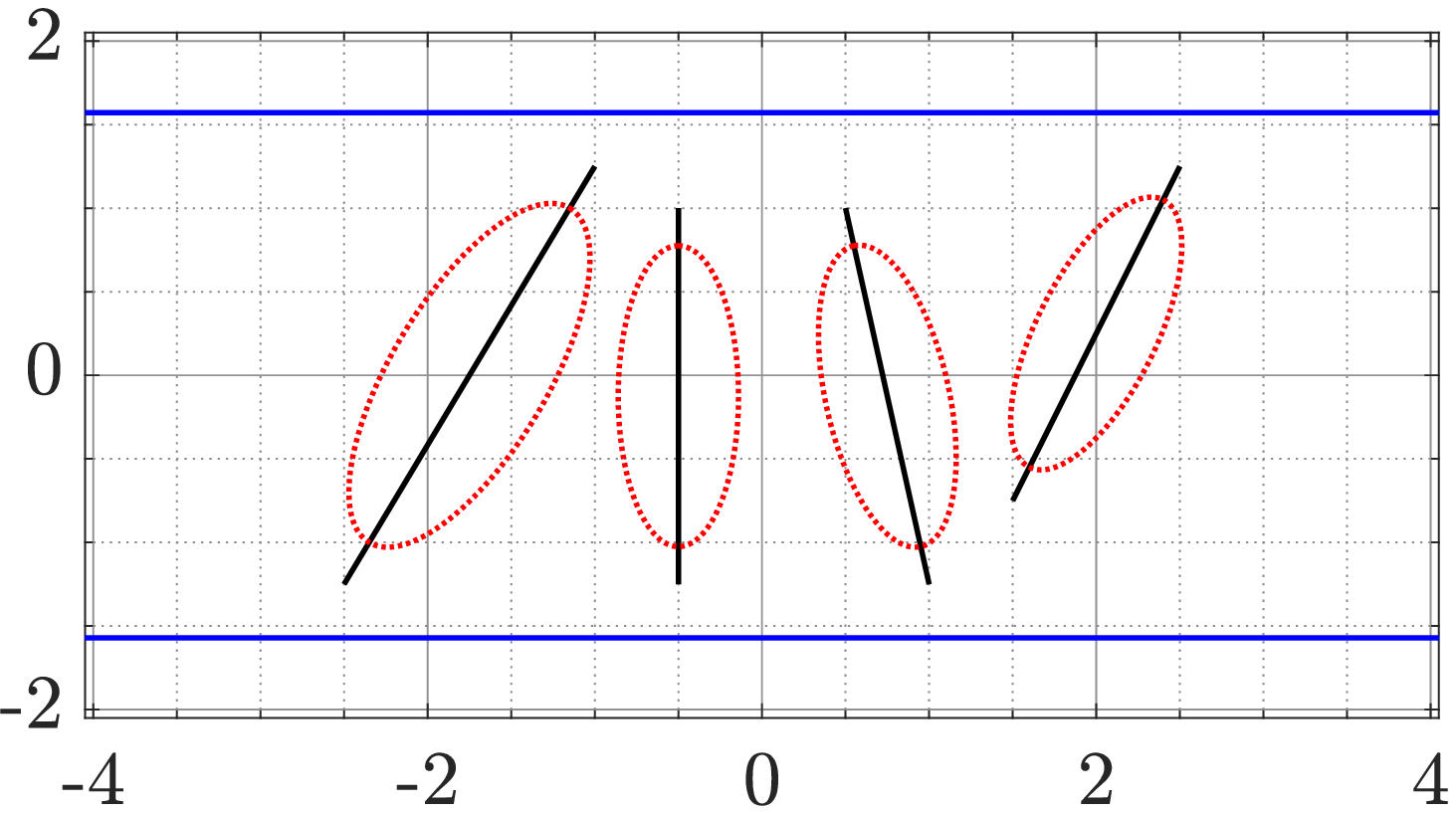}}
\hfill
\scalebox{0.3}{\includegraphics[trim=0 0    0 0,clip]{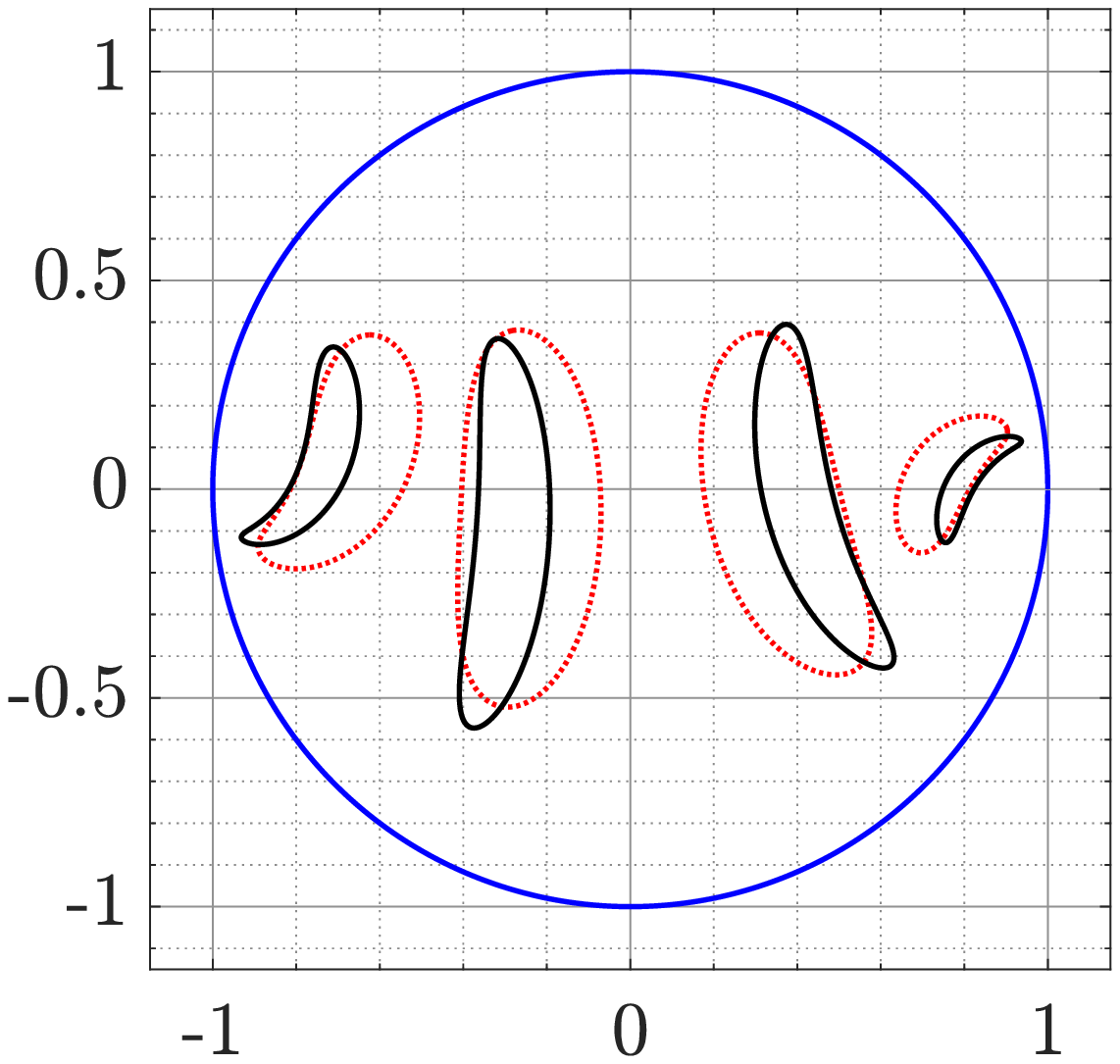}}
}
\caption{On the left, the domain $\Omega$ with four rectilinear slits. In the center, the initial intermediate domain $\hat\Omega^0$ where thin ellipses are chosen on the given slits with $r=0.4$. On the right, the domain interior to the unit circle and exterior to the dot-red curves is the initial domain $G^0$. The computed preimage domain $G$ is the domain interior to the unit circle and exterior to the solid-black curves.}
\label{fig:strip}
\end{figure}

For $k=0,1,2,3,\ldots$, where $k$ is the iteration number, we let $\hat\Omega^k$ to be the multiply connected domain obtained by removing the $m$ ellipses $\hat L_j^k$ parametrized by
\[
\hat\eta_j^k(t)= z_j^k+0.5a_j^ke^{\i\theta_j}(\cos t-\i r\sin t), \quad t\in J_j, \quad j=1,\ldots,m,
\]
from the infinite strip $S=\{\xi\in\C\,:\, \left|\Im \xi\right|<\pi/2\}$. Here $r$, $0<r\le1$, is the ratio of the lengths of the major to the minor axes of the ellipse $\hat L_j^k$. 
We denote by $G^k$ the image of the multiply connected domain $\hat\Omega^k$ under the conformal mapping $w=\Psi^{-1}(\xi)$ where 
\[
w=\Psi^{-1}(\xi)=\frac{\exp(\xi)-1}{\exp(\xi)+1}=\tanh(\xi/2).
\]
Hence, the preimage domain $G^k$ (see Figure~\ref{fig:strip} (right) for an example) is the bounded multiply connected domain interior to the unit circle parametrized by
\[
\eta^k_0(t)=e^{\i t}, \quad  t\in J_0,
\]
and exterior to the $m$ quasi-ellipses $\Gamma_1,\ldots,\Gamma_m$ parametrized by
\[
\eta^k_j(t)=\Psi^{-1}\left(\hat\eta_j^k(t)\right), \quad t\in J_j, \quad j=1,\ldots,m.
\]
Thus, determining the the preimage domain $G^k$ requires computing the parameters $z_j^k$ and $a_j^k$ of the ellipses $\hat L_j^k$, $j=1,\ldots,m$, which will be accomplished using the following iterative method. 
We point out that the initial domain $\hat\Omega^0$ is obtained from the given domain $\Omega$ by replacing each slit $L_j$ by a thin ellipse $\hat L_j^0$ whose major axis is on the slit $L_j$ (see Figure~\ref{fig:strip} (center)).\\

\noindent{\bf Initialization:}\\
Set
\[
z_j^0=c_j, \quad a_j^0=(1-0.5r)\ell_j, \quad j=1,\ldots,m.
\]
\noindent{\bf Iterations:} \\
For $k=1,2,3,\ldots$,
\begin{itemize}
	\item The method in Theorem~\ref{thm:cm-strip} is used to compute the conformal mapping from the domain $G^{k-1}$ onto the domain $\Omega^k$ (the infinite strip $|\Im z|<\pi/2$ with $m$ slits $L_1^k,\ldots,L_m^k$ such that the angle between the slit $L_j^k$ and the positive real axis is $\theta_j$, $j=1,\ldots,m$).
	\item If $c_j^k$ is the center of the slit $L_j^k$ and $\ell_j^k$ is its length, then the parameters $z_j^k$ and $a_j^k$ are updated through
\begin{eqnarray*}
     z_j^{k} &=& z_j^{k-1}-(c_j^{k}-c_j), \\
     a_j^{k} &=& a_j^{k-1}-(1-0.5r)(\ell_j^k -\ell_j),
\end{eqnarray*}
for $j=1,\ldots,m$.
  \item Stop the iterations if
	\begin{equation}\label{eq:Ek}
	E_k=\frac{1}{2m}\sum_{j=1}^{m}\left(|c_j^{k} -c_j|+|\ell_j^k -\ell_j|\right)<\varepsilon \quad{\rm or}\quad k>{\tt Max}
	\end{equation}
	where $\varepsilon$ is a fixed tolerance and ${\tt Max}$ is the maximum number of iterations to be not exceeded. In our numerical experiments, we set $\varepsilon=10^{-14}$ and ${\tt Max}=100$.
\end{itemize}

The iterations above produce a sequence of multiply connected preimage domains $G^0$, $G^1$, $G^2$, $G^3$, $\ldots$, which numerically converges to the required domain $G$. This method also produces a conformal map $z=\Phi(w)$ from the computed preimage domain $G$ onto the given domain $\Omega$. Similar iterative procedures have been experimentally studied for other canonical domains in~\cite{NG18,N19} where the numerical examples demonstrate the fast convergence even for domains with high connectivity. 

It is understandable that at every iteration of the above method one needs to solve the integral equation~\eqref{eq:ie} and to compute the piecewise constant function $h$ in~\eqref{eq:h}. This can be done with the fast method presented in~\cite{Nas-ETNA} by applying the MATLAB function \verb|fbie| in which we discretize~\eqref{eq:ie} by the Nystr\"om method using the trapezoidal rule with $n$ equidistant nodes in each sub-interval $J_j$, $j=0,\ldots,m$. This produces a $(m+1)n\times(m+1)n$ linear system, which in turn is solved by the generalized minimal residual method. 
More precisely, we employ the MATLAB function \verb|gmres| without restart where the tolerance and maximum number of iterations are chosen to be $10^{-14}$ and $100$, respectively. 
The matrix-vector product in \verb|gmres| is computed through the Fast Multipole Method.
by calling \verb|zfmm2dpart| from the MATLAB toolbox~\cite{Gre-Gim12} with a  tolerance of $0.5\times10^{-15}$.
The complexity of each iteration of the presented iterative method is $O((m+1)n\log n)$.
We refer to ~\cite{Nas-ETNA} for further details.

For the parameter $r$ in the above iterative method, in general, the value of $0<r\le 1$ is chosen such that the inner boundary components are not overlapping, i.e., we need to choose a small value for $r$ if the slits are close to each others. On one hand, the geometry of the preimage domain $G$ will be simpler if we choose $r=1$ or close to $1$ (when it is possible). On the other hand, the rate of convergence of the iterative method depends on the value of $r$ where the method converges faster for small $r$ (see Figure~\ref{fig:ex-m22-err} (left) below for the example considered in Section~\ref{ex:cm}). However, for small values of $r$, the method requires more GMRES iterations (see Figure~\ref{fig:ex-m22-err} (right)). Furthermore, for small values of $r$, the inner boundary components will be thin and hence one needs to consider a larger value of $n$ to achieve a satisfactory accuracy. The preimage domain in Figure~\ref{fig:ex-m22} is computed with $r=0.2$.

\subsection{Computing the conformal mapping from $\Omega$ onto $G$}

The iterative method presented in the preceding subsection allows us to compute the preimage domain $G$ and the conformal mapping from $z=\Phi(w)$ from $G$ onto $\Omega$. More precisely, it provides us with a parametrization $\eta(t)$, $t\in J$, of the boundary $\partial G$. As described in Theorem~\ref{thm:cm-strip}, solving~\eqref{eq:ie} and computing $h$ in~\eqref{eq:h} permit to determine the boundary values of the auxiliary analytic function $f$.
Consequently, the boundary values $\Phi(\eta(t))$ of $\Phi$ can be computed by~\eqref{eq:Phi-f}. 
To determine the values $\Phi(w)$ for $w\in G$, we first compute the values of $f(w)$ using the Cauchy integral formula. In our numerical computations, the values $f(w)$ 
for $w\in G$ are computed accurately using the MATLAB function \verb|fcau| presented in~\cite{Nas-ETNA}. 

\begin{figure}[ht] %
\centerline{
\scalebox{0.3}{\includegraphics[trim=0 0    0 0,clip]{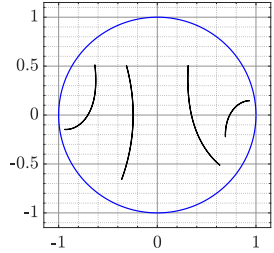}}
}
\caption{The bounded domain $\tilde\Omega$ for the domain $\Omega$ shown in Figure~\ref{fig:strip}(left).}
\label{fig:t-omega}
\end{figure}

Since $\Phi$ maps $\partial G$ onto $\partial\Omega$, we see that
\begin{equation}\label{eq:zet}
\zeta(t) = \Phi(\eta(t)) = -\i f(\i)+\eta(t) f(\eta(t))+\Psi(\eta(t)), \quad t\in J,
\end{equation}
is a parametrization of $\partial\Omega$, and $\Phi^{-1}(\zeta(t))=\eta(t)$. The boundary of the domain $\Omega$ passes through the point at infinity, which means the Cauchy integral formula cannot be employed to directly compute $\Phi^{-1}(z)$ for $z\in\Omega$. Nevertheless, we write $\Phi^{-1}(z)$ as
\begin{equation}\label{eq:Phi-1-g}
\Phi^{-1}(z)=(\Phi^{-1}\circ\Psi\circ\Psi^{-1})(z) = g(\Psi^{-1}(z))
\end{equation}
where $g=\Phi^{-1}\circ\Psi$ and $\Psi$ is given by~\eqref{eq:Psi}. The mapping function $\Psi^{-1}$ maps the unbounded domain $\Omega$ onto the bounded domain $\tilde\Omega=\Psi^{-1}(\Omega)$ (see Figure~\ref{fig:t-omega}). The domain $\tilde\Omega$ is interior to the unit circle and exterior to $m$ slits (which are not rectilinear) and the boundary of $\tilde\Omega$ is parametrized by 
\[
\tilde\zeta(t) = \Psi^{-1}(\zeta(t)), \quad t\in J.
\]
The boundary values of $g$ are then given by
\[
g(\tilde\zeta(t))=\Phi^{-1}(\Psi(\tilde\zeta(t)))=\Phi^{-1}(\zeta(t))=\eta(t), \quad t\in J.
\]
For $z\in\Omega$ we have $\tilde z=\Psi^{-1}(z)\in\tilde\Omega$, and the values of $g(\tilde z)$ can now be computed using the Cauchy integral formula.
Consequently, the values of $\Phi^{-1}(z)$ can be computed through~\eqref{eq:Phi-1-g}, that is,
\[
\Phi^{-1}(z)= g(\Psi^{-1}(z)) 
= \frac{1}{2\pi\i} \int_{\partial\tilde\Omega} \frac{g(\tilde\zeta)}{\tilde\zeta-\Psi^{-1}(z)}d\tilde\zeta
= \frac{1}{2\pi\i} \int_{J} \frac{\eta(t)}{\tilde\zeta(t)-\Psi^{-1}(z)}\tilde\zeta'(t)dt.
\]

In the numerical computations presented below, the values of $g(\tilde z)$ are computed by the MATLAB function \verb|fcau|. However, to use the function \verb|fcau|, we need to compute $\tilde\zeta'(t)$ for $t\in J$ in advance. This derivative is computed by first approximating the real and imaginary parts of $\tilde\zeta(t)$ on each sub-interval $J_j$, $j=0,\ldots,m$, with trigonometric interpolating polynomials, and then differentiating these polynomials. Note that the interpolating polynomials can be computed by applying the fast Fourier transform~\cite{Weg05}.

All computations in this paper are performed in MATLAB R2017a on an ASUS Laptop with Intel Core i7-8750H CPU @ 2.20GHz, 2208 Mhz, 6 Cores, 12 Logical Processors and 16 GB RAM. 

\subsection{A numerical example}\label{ex:cm}

\begin{figure}[ht] %
	\centerline{
		\scalebox{0.35}{\includegraphics[trim=0 0 0 0,clip]{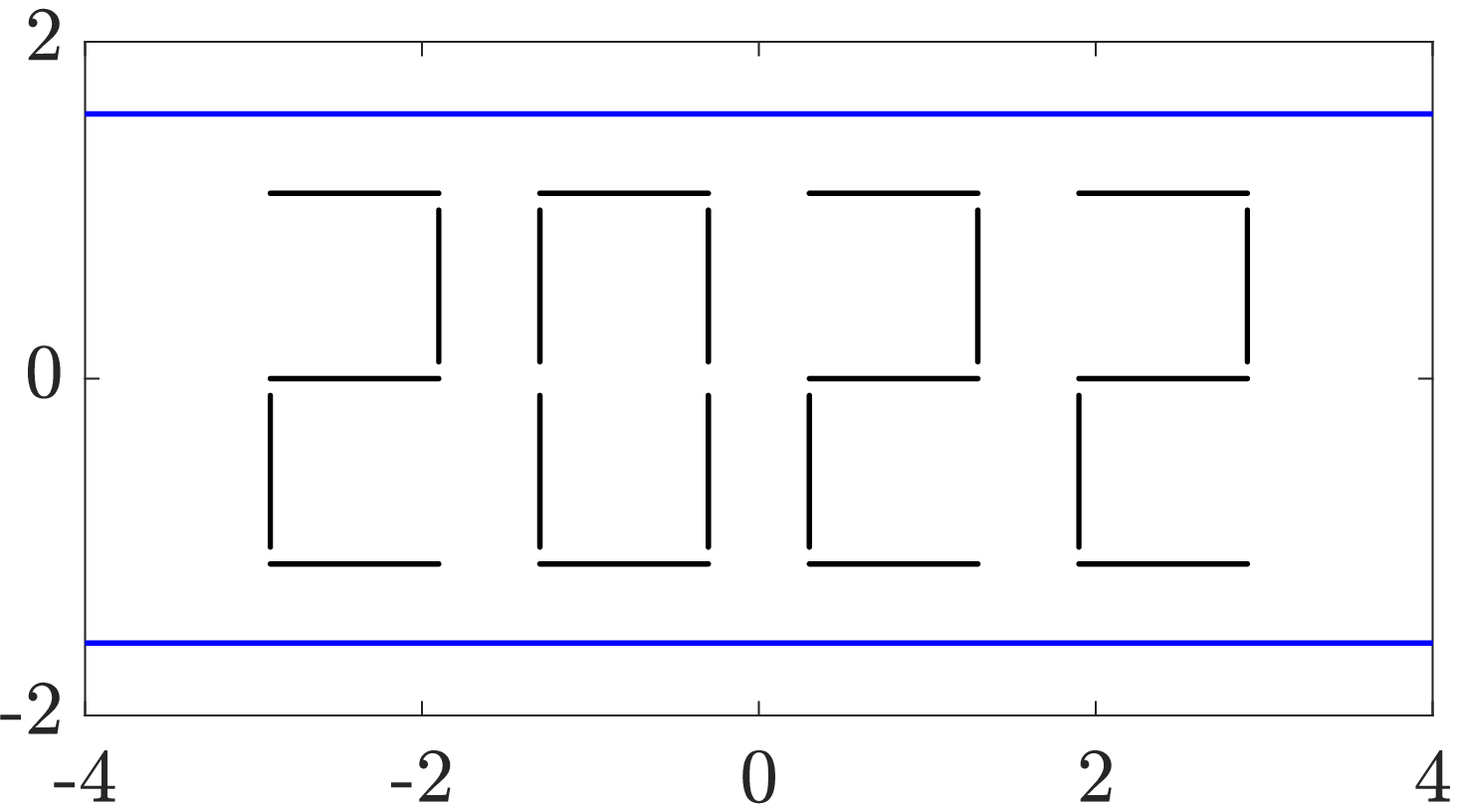}}
		\hfill
		\scalebox{0.35}{\includegraphics[trim=0 0 0 0,clip]{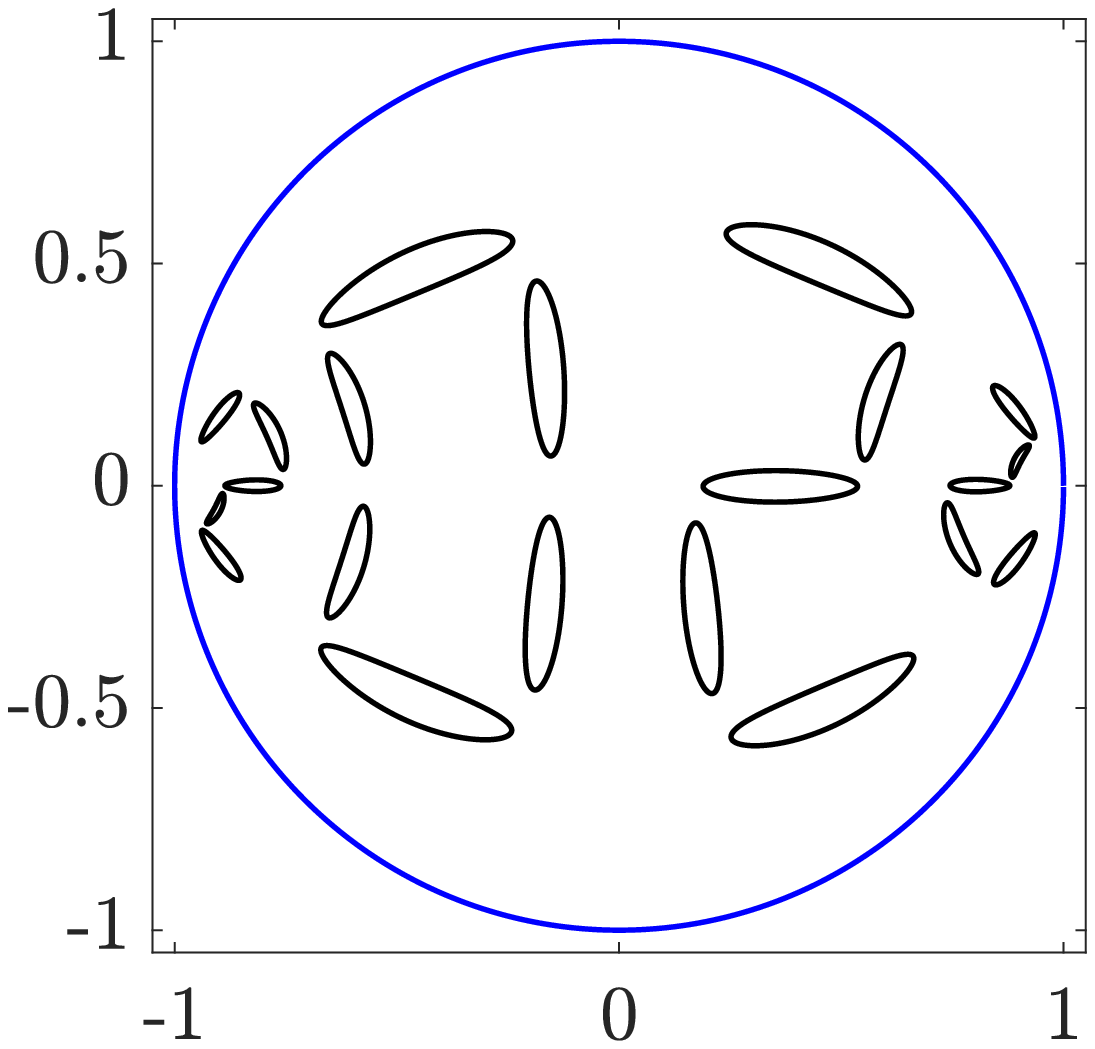}}
	}
	\caption{The domain $\Omega$ with $21$ rectilinear slits (left) and the preimage domain $G$ (right).}
	\label{fig:ex-m22}
\end{figure}

\begin{figure}[ht] %
	\centerline{
		\scalebox{0.35}{\includegraphics[trim=0 0.0cm  0 0,clip]{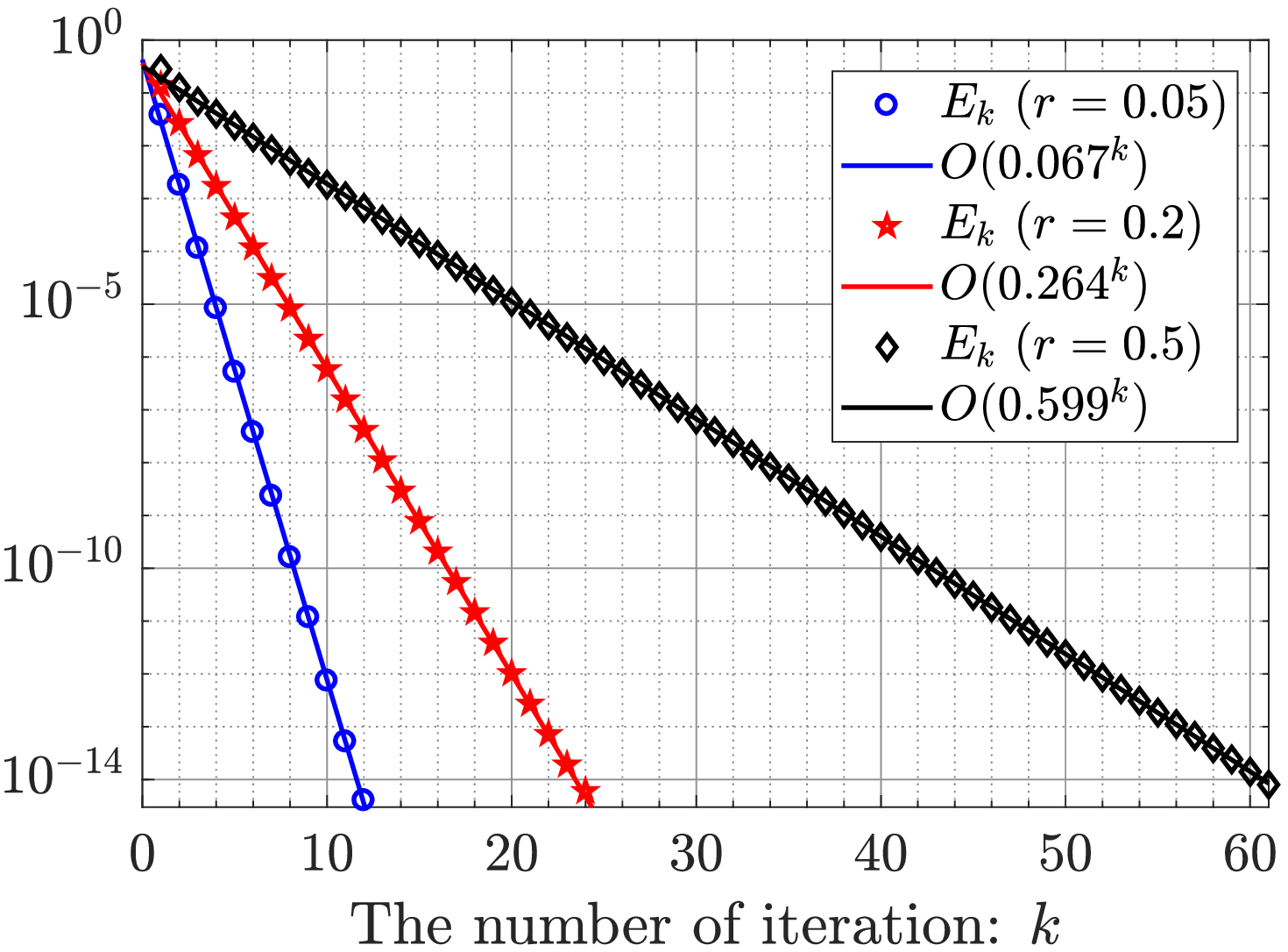}}
		\hfill
		\scalebox{0.35}{\includegraphics[trim=0 0.0cm  0 0,clip]{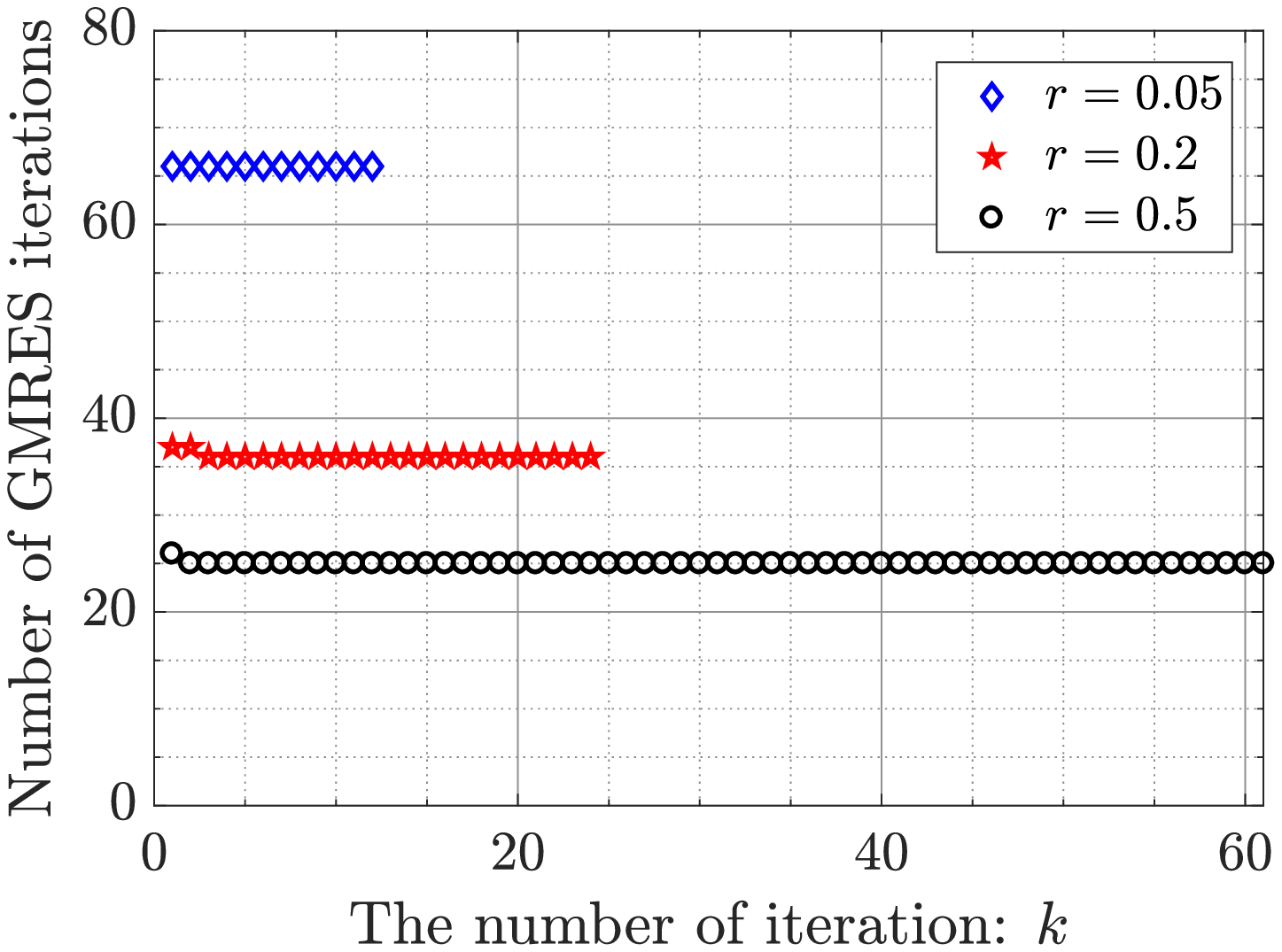}}
	}
	\caption{The error $E_k$ for the domain in Figure~\ref{fig:ex-m22} (left) for $r=0.1$, $r=0.2$, and $r=0.3$.}
	\label{fig:ex-m22-err}
\end{figure}

We present a numerical example that illustrates how the iterative method introduced above can be applied to determine the preimage $G$ of an infinite strip with $m=21$ rectilinear slits as shown in Figure~\ref{fig:ex-m22} (left). The method is applied with $n=2^{11}$ and $r=0.2$ and the obtained preimage domain is displayed in Figure~\ref{fig:ex-m22} (right). We also compute the error $E_k$ defined by~\eqref{eq:Ek} and the number of GMRES iterations using $n=2^{11}$ for $r=0.05$, $r=0.2$, and $r=0.5$. The computed error and the number of GMRES iterations vs. the number of iteration $k$ is given in Figure~\ref{fig:ex-m22-err}. The total CPU time required by the iterative method is $68.91$ sec for $r=0.05$, $70.15$ sec for $r=0.2$, and $127.72$ for $r=0.5$. 

\section{Capacity of generalized condensers}
\label{sec:strip-g-cap}

\subsection{Numerical computation of capacity of generalized condensers}

\begin{figure}[ht] %
\centerline{
\scalebox{0.4}{\includegraphics[trim=0 0 0 0,clip]{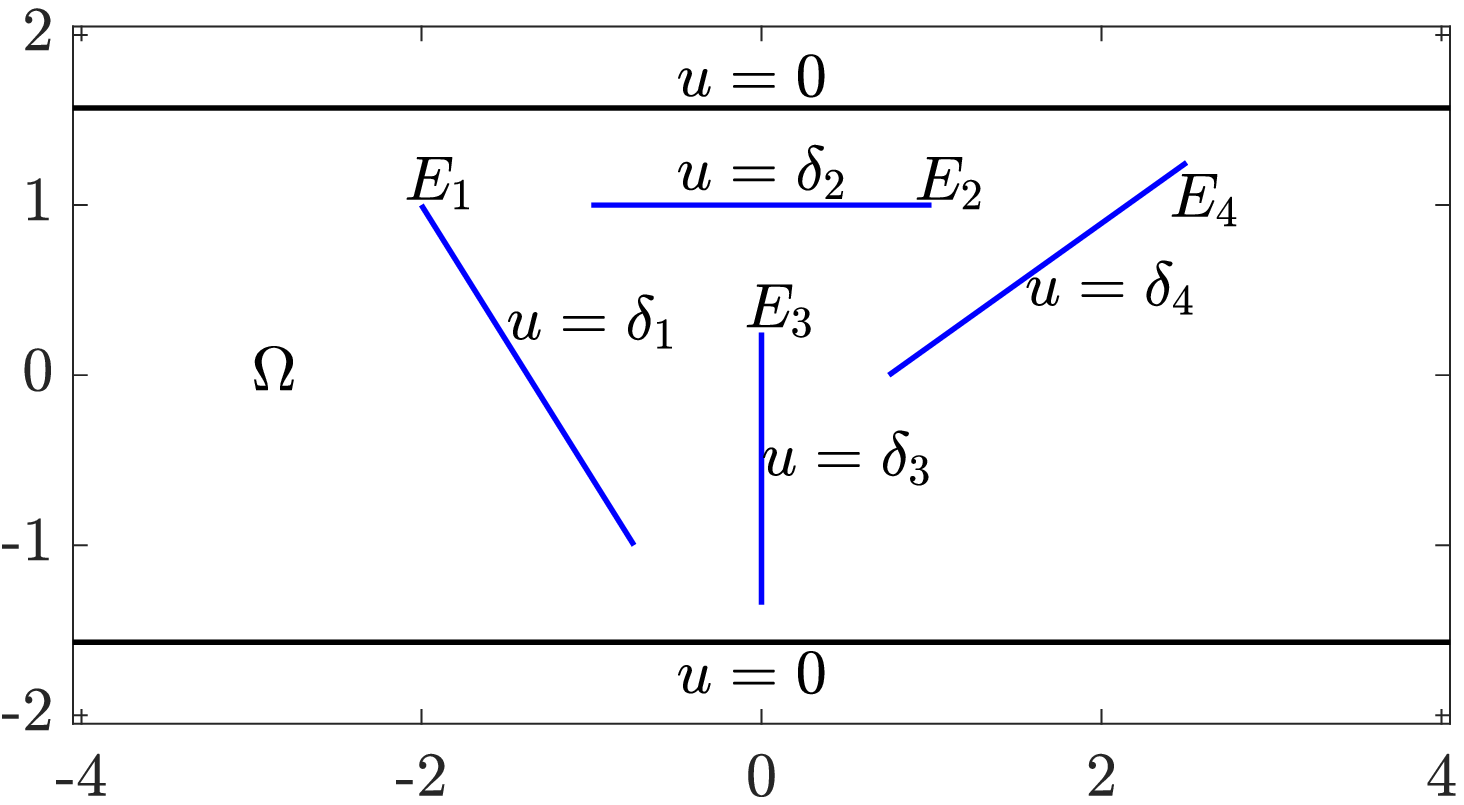}}
\hfill
\scalebox{0.4}{\includegraphics[trim=0 0 0 0,clip]{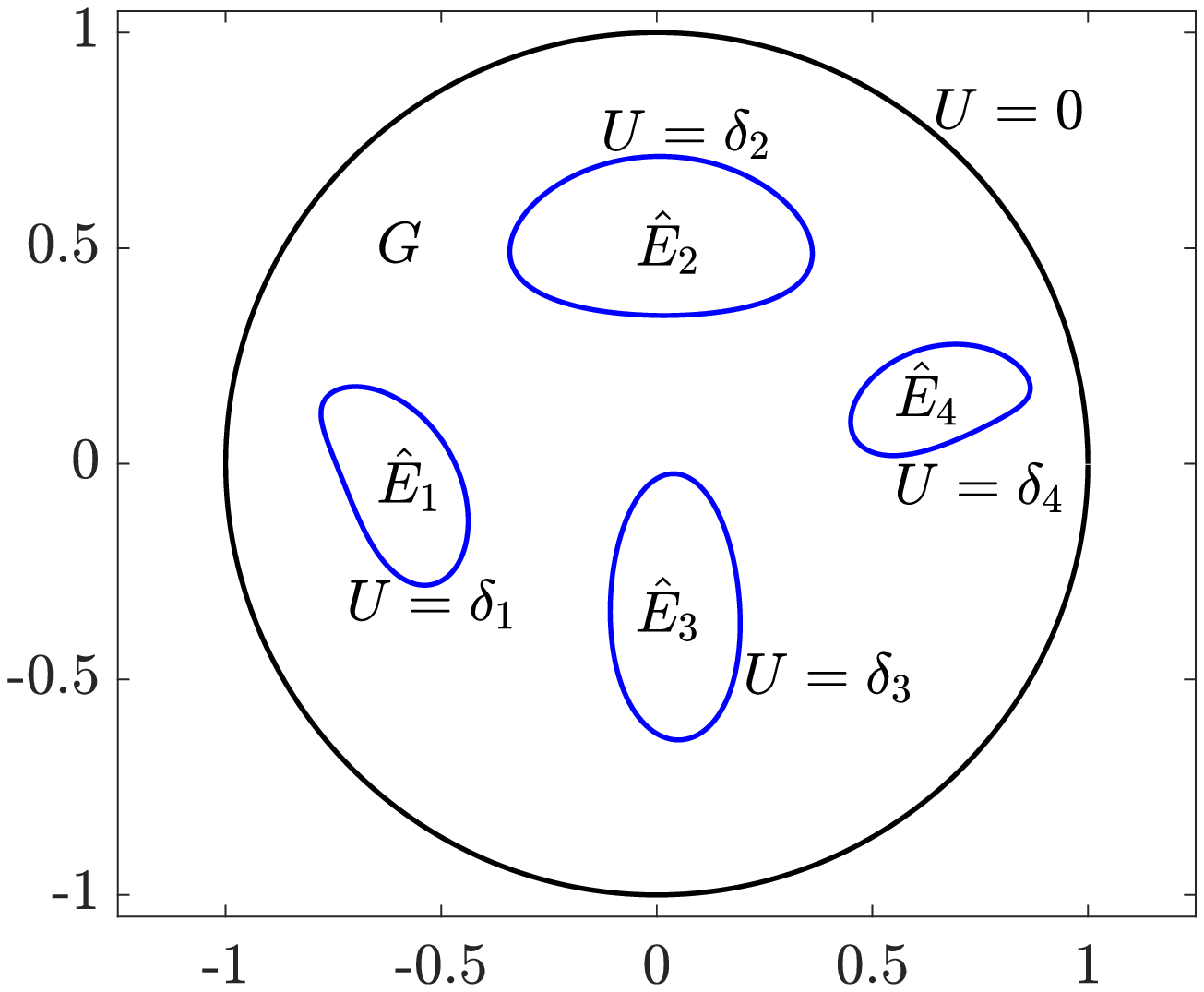}}
}
\caption{The domain $\Omega$ (left) and the domain $G$ (right) for $m=4$.}
\label{fig:domaingcS}
\end{figure}

We aim to compute the capacity of generalized condensers of the form $C=(S,E,\delta)$ where $S$ is the infinite strip 
\[
S=\{z\,:\, \left|\Im z\right|<\pi/2\},
\]
$E=\{E_j\}_{j=1}^m$ is a collection of $m$ nonempty closed pairwise disjoint segments $E_j=[a_j,b_j]$ with complex numbers $a_j,b_j\in S$, and $\delta=\{\delta_j\}_{j=1}^{m}$ is a collection of real numbers. We assume that $m>1$ and $\delta$ contains at least two different numbers. The domain $\Omega=S \setminus E =S \setminus \cup_{j=1}^{m}[a_j,b_j]$ is known as the field of the condenser $C$, the sets $E_k$ as the plates, and the numbers $\delta_k$ as the levels of the potential of the plates $E_j$, $j=1,\ldots,m$~\cite[p.~12]{du}  (see Figure~\ref{fig:domaingcS} (left) for $m=4$).

The conformal capacity of the generalized condenser $C$, $\capa(C)$, is given by the Dirichlet integral
\begin{equation}\label{eq:cap}
\capa(C)=\iint\limits_{\Omega}|\nabla u|^2dxdy
\end{equation}
where $u$ is the potential function of the condenser $C$, i.e., 
the function $u$ is the unique solution of the Dirichlet BVP~\cite[p.~13, p.~305]{du}
\begin{subequations}\label{eq:u-Dir}
\begin{eqnarray}
\label{eq:u-Dir-1}
\Delta u&=& 0, \quad\; {\rm on}\;\;\Omega, \\
\label{eq:u-Dir-2}
u      &=& 0, \quad\; {\rm on}\;\;\partial S, \\
\label{eq:u-Dir-3}
u      &=& \delta_j, \quad {\rm on}\;\;E_j, \quad j=1,\ldots,m.
\end{eqnarray}
\end{subequations}

If all values of $\delta_k$ are equal to $1$ or if $m=1$ (i.e., $E$ consists of only one segment), then $C=(S,E,\delta)$ is a classical condenser and it is shortly denoted by $C=(S,E)$. In this case, the capacity $\capa(C)$ may simply be written as $\capa(S,E)$, see Section~\ref{section1}. 

One important property of the capacity of generalized condensers is its invariance  under conformal mappings, which implies that domains having complex geometry like $G$  can be treated with the aid of conformal mappings~\cite{ah,avv,bbgghv,crow-20,dnv,dek,dt,du,et,garmar,hkv,LV,Nvm,Nv,ps,sl,sol,Vas02,vu88,Weg05,wen}.

The capacity $\capa(C) = \capa(S,E,\delta)$ will be computed in two steps:\\
\noindent{\bf Step 1.}  For the given domain $\Omega=S \setminus \cup_{j=1}^{m}[a_j,b_j]$, the iterative method proposed in Section~\ref{sec:strip} will be used to find a preimage domain $G$ interior to the unit circle $\Gamma_0$ and exterior to $m$ ellipses $\Gamma_j$, $j=1,\ldots,m$. Let $D_j$ be the simply connected domain interior to $\Gamma_j$ for $j=1,\ldots,m$, let $\hat E=\cup_{j=1}^{m}\overline{D}_j$, and let $\hat C=(\D,\hat E,\delta)$ (see Figure~\ref{fig:domaingcS} (right)). Then $\capa(C)=\capa(\hat C)$.\\
\noindent{\bf Step 2.}  The capacity $\capa(\hat C)$ will be evaluated by the method that we will describe in the sequel.

We point out that we use the same number $n$ of discretization points for both steps. 
Now, let us present the method needed in Step 2 to compute the capacity of the generalized condenser $\hat C=(\D,\hat E,\delta)$ where $\D$ is the unit disk, $\hat E=\{\hat E_j\}_{j=1}^m$ is a family of $m$ nonempty closed and pairwise disjoint sets $\hat E_j=D_j\cup\Gamma_j\subset \D$, and $\delta=\{\delta_j\}_{j=1}^{m}$ is a collection of real numbers. Here, $\Gamma_j=\partial \hat E_j$ is a smooth Jordan curve for $j=1,\ldots,m$. Hence, $G=\D \setminus \hat E$ is a bounded multiply connected domain of connectivity $m+1$, and $\Gamma=\partial G=\cup_{j=0}^{m}\Gamma_j$ where $\Gamma_0$ is the unit circle and $\Gamma_1,\ldots,\Gamma_m$ are ellipses enclosed in $\Gamma_0$ (see Figure~\ref{fig:domaingcS} (right) for $m=4$). 

The above generalized condenser $\hat C=(\D,\hat E,\delta)$ is a special type of the generalized condenser considered in~\cite[p.~12]{du} and~\cite{Nvm}. Hence, the numerical method presented in~\cite{Nvm} can be used to compute the capacity of the above generalized condenser. For the  convenience of the reader, we review this method below.

The conformal capacity of the generalized condenser $\hat C$ is
\begin{equation}\label{eq:capU}
\capa(\hat C)=\iint\limits_{G}|\nabla U|^2dxdy
\end{equation}
where the potential function $U$ is now the unique solution of the Dirichlet BVP
\begin{subequations}\label{eq:U-Dir}
\begin{eqnarray}
\label{eq:U-Dir-1}
\Delta U&=& 0, \quad\; {\rm on}\;\; G, \\
\label{eq:U-Dir-2}
U      &=& 0, \quad\; {\rm on}\;\;\Gamma_0, \\
\label{eq:U-Dir-3}
U      &=& \delta_j, \quad {\rm on}\;\;\Gamma_j, \quad j=1,\ldots,m.
\end{eqnarray}
\end{subequations}
The harmonic function $U$ is the real part of an analytic function $F$ in $G$ which is not necessarily single-valued. If we take a given point $\alpha_k\in G_k$ for each $k=1,\ldots,m$, then $F$ can be written as~\cite{Gak,garmar,Mik64,Mus}
\begin{equation}\label{eq:F-u}
F(z)=g(z)-\sum_{k=1}^{m} a_k\log(z-\alpha_k)
\end{equation}
where $g$ is a single-valued analytic function in $G$ and $a_1,\ldots,a_{m}$ are undetermined real constants satisfying~\cite[\S31]{Mik64}
\begin{equation}\label{eq:ak}
a_k=\frac{1}{2\pi}\int_{\Gamma_k}\frac{\partial U}{\partial{\bf n}}ds, \quad k=1,\ldots,m.
\end{equation}
Using Green's formula~\cite[p.~4]{du}, and in view of~\eqref{eq:U-Dir-2}--\eqref{eq:U-Dir-3} and~\eqref{eq:ak}, equation~\eqref{eq:capU} can be written as
\begin{equation}\label{eq:cap-n}
\capa(\hat C)=\int_{\Gamma}u\frac{\partial u}{\partial{\bf n}}ds
=\sum_{k=1}^{m}\delta_k\int_{\Gamma_k}\frac{\partial u}{\partial{\bf n}}ds=2\pi\sum_{k=1}^{m}\delta_k a_k.
\end{equation}
Consequently, by the conformal invariance of the capacity, we have
\begin{equation}\label{eq:cap-a}
\capa(C)=2\pi\sum_{k=1}^{m}\delta_k a_k.
\end{equation}

The problem~\eqref{eq:U-Dir} above is a particular case of the problem considered in~\cite[Eq.~(4)]{Nvm}. Let $\eta(t)$, $t\in J$, be a parametrization of the boundary $\Gamma=\partial G$ and let $A(t)$ be defined~\eqref{eq:A}. Then, the constants $a_1,\ldots,a_m$ in~\eqref{eq:cap-a} will be computed as in the following theorem from~\cite[Theorem 4]{Nvm}. 

\begin{theorem}\label{thm:method}
For each $k=1,\ldots,m$, let the function $\gamma_k$ be defined by
\begin{equation}\label{eq:gam-k}
\gamma_k(t)=\log|\eta(t)-\alpha_k|,
\end{equation}
let $\rho_k$ be the unique solution of the integral equation
\[
(\bI-\bN)\rho_k=-\bM\gamma_k,
\]
and let the piecewise constant function $h_k=(h_{0,k},h_{1,k},\ldots,h_{m,k})$ be given by
\[
h_k=[\bM\rho_k-(\bI-\bN)\gamma_k]/2.
\]
Then, the $m+1$ real constants $a_1,\ldots,a_{m},c$ are the unique solution of the linear system
\begin{equation}\label{eq:sys-method}
\left[\begin{array}{ccccc}
h_{0,1}    &h_{0,2}    &\cdots &h_{0,m}      &1       \\
h_{1,1}    &h_{1,2}    &\cdots &h_{1,m}      &1       \\
\vdots     &\vdots     &\ddots &\vdots       &\vdots  \\
h_{m,1}    &h_{m,2}    &\cdots &h_{m,m}      &1       \\
\end{array}\right]
\left[\begin{array}{c}
a_1    \\a_2    \\ \vdots \\ a_{m} \\  c 
\end{array}\right]
= \left[\begin{array}{c}
0 \\  1 \\  \vdots \\ 1  
\end{array}\right].
\end{equation}
\end{theorem}

The linear system~\eqref{eq:sys-method} is usually small and can be solved using Gaussian elimination. Henceforth, the capacity can be computed by~\eqref{eq:cap-a}. Note that, in our computation in this paper, we do not need the value of the real constant $c$ in~\eqref{eq:sys-method}.

We now present several examples to illustrate how the above proposed method can be applied to approximate the capacity $\capa(S,E,\delta)$. 

\subsection{A strip with one rectilinear slit} 

\begin{example}\label{ex:msi-si}
$E=[-s\i,s\i]$ where $0<s<\pi/2$ is a real number.
\end{example}

The exact value of the capacity of the condenser $(S,E)$ can be given in terms of special functions 
for this example. For $\Psi$ defined by~\eqref{eq:Psi}, the mapping function
\[
z\mapsto-\i\Psi^{-1}(z)=-\i\tanh(z/2)
\]
maps conformally the domain $S\backslash[-s\i,s\i]$ onto the domain $\D\backslash[-\tan(s/2),\tan(s/2)]$. Then, by the M{\"o}bius transformation
\[
z\mapsto \frac{z+\tan(s/2)}{1+z\,\tan(s/2)}
\]
the domain $\D\setminus[-\tan(s/2),\tan(s/2)]$ is mapped onto the domain $\D\setminus[0,\sin(s)]$. 
Owing to the conformal invariance of the capacity, we have
\[
\capa(S,[-s\i,s\i]) = \capa(\D,[0,\sin(s)]).
\]
This yields (see \cite{LV}, \cite[Thm 8.6(1)]{vu88}),
\begin{equation}\label{eq:strip-exact-i}
\capa(S,[-s\i,s\i])=\frac{2\pi}{\mu(\sin s)}
\end{equation}
where 
\[
\mu(r)=\frac{\pi}{2}\frac{\K'(r)}{\K(r)}, \quad
\K(r)=\int^1_0 \frac{dx}{\sqrt{(1-x^2)(1-r^2x^2)}},
\quad \K'(r)=\K\left(\sqrt{1-r^2}\right).
\]
Here, $\K(r)$ and $\K'(r)$ are the elliptic integrals of the first kind, and $\mu:(0,1)\to(0,\infty)$ is a decreasing homeomorphism. For the numerical computation of the values of $\mu(r)$, we use the method described in~\cite{Nv}.

\begin{figure}[ht] %
\centerline{
\scalebox{0.4}{\includegraphics[trim=0 0.0cm 0 0cm,clip]{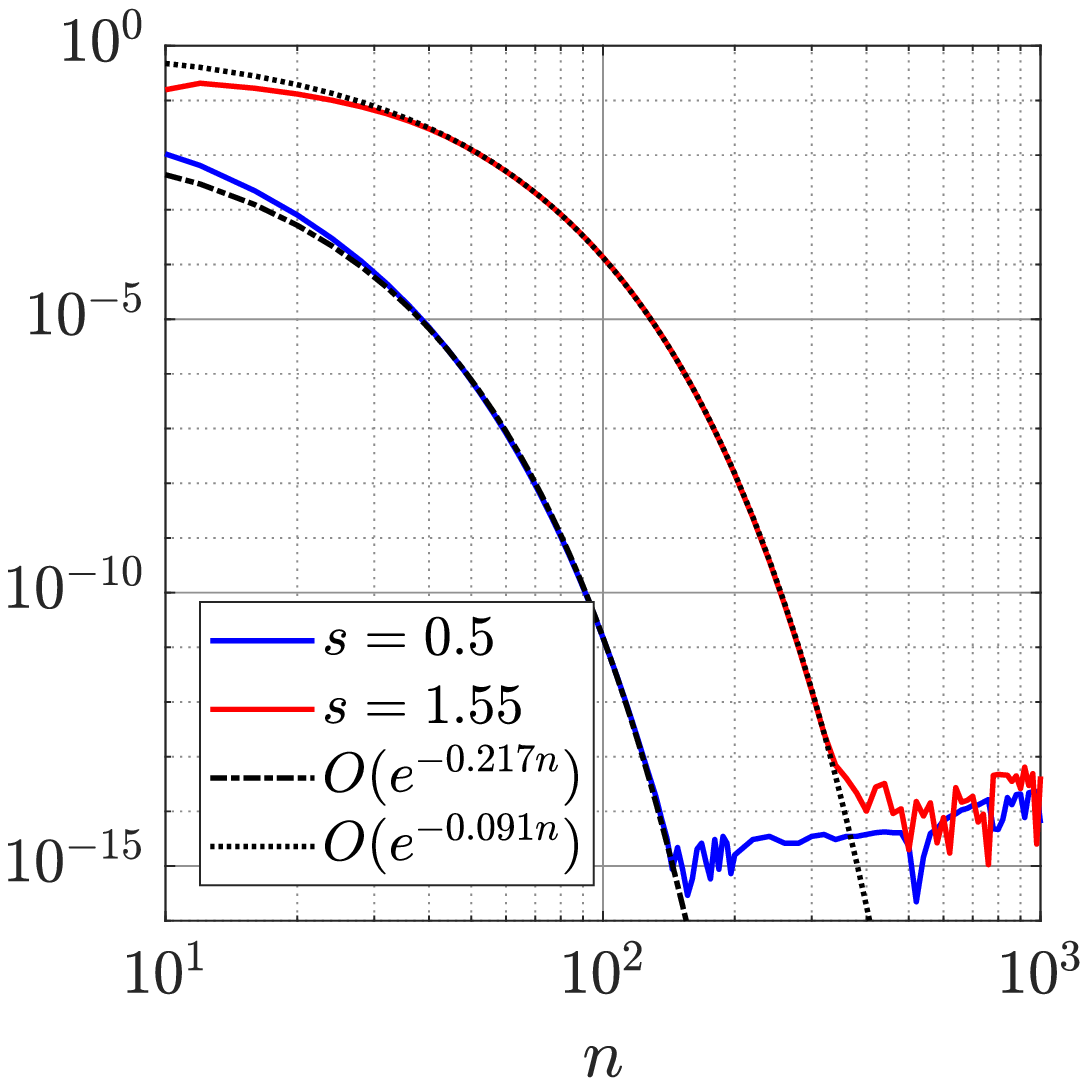}}
\hfill
\scalebox{0.4}{\includegraphics[trim=0 0.0cm 0 0cm,clip]{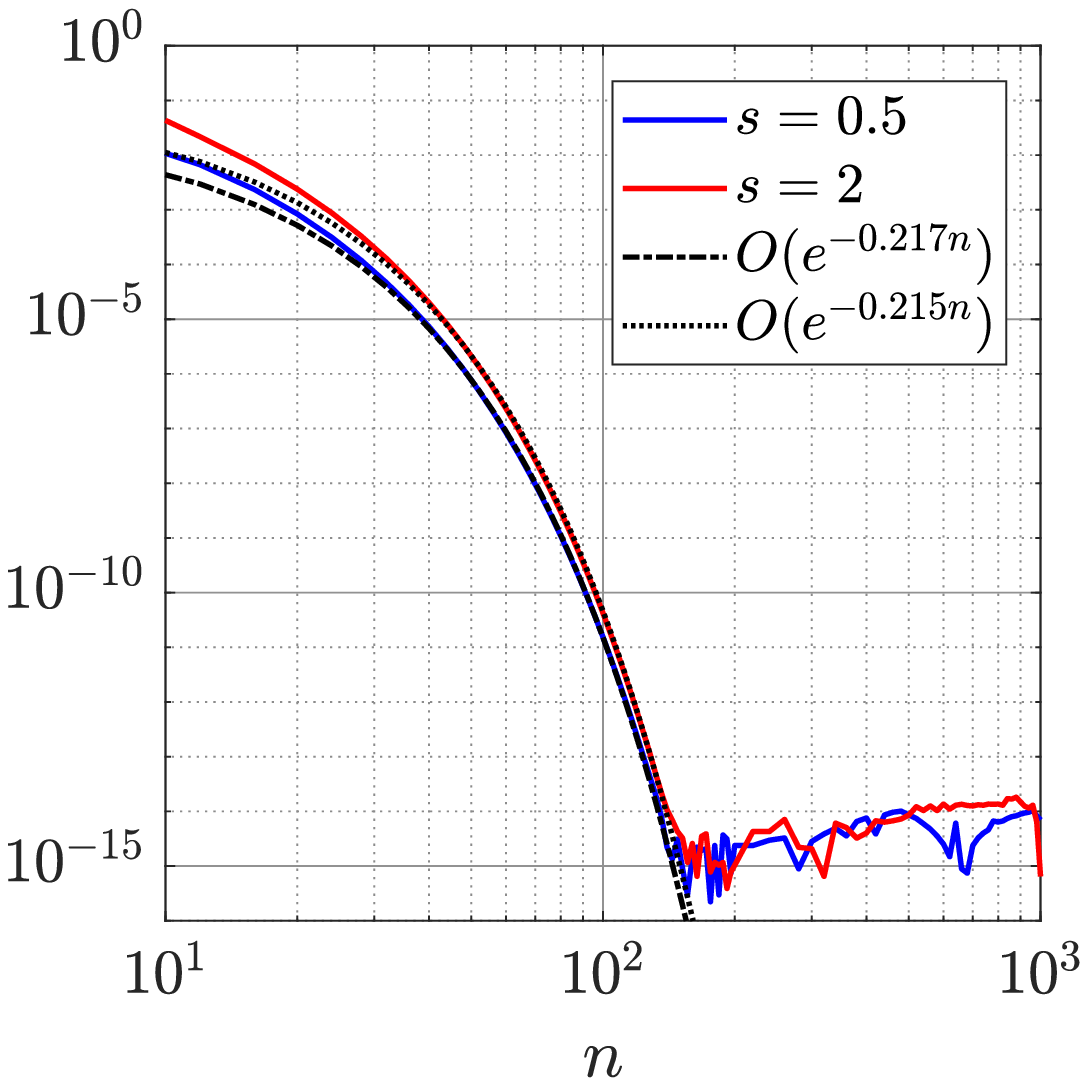}}
}
\caption{On the left, the relative error in the approximated values for Example~\ref{ex:msi-si} for $s=0.5$ and $s=1.55$. On the right, the relative error in the approximated values for Example~\ref{ex:ms-s} for $s=0.5$ and $s=2$.}
\label{fig:sing-msi-err}
\end{figure}

We set $n=2^{10}$ and $r=0.2$ to compute $\capa(S,E)$ for several values of $s$ in $(0,\pi/2)$, and present the numerical results in Table~\ref{tab:strip-seg-err-i}. 
The relative error in the approximate values for $s=0.5$ and $s=1.55$ vs. $n$, the number of discretization points on each boundary component, are depicted in Figure~\ref{fig:sing-msi-err} (left). The error is $O(e^{-0.217n})$ for $s=0.5$ and $O(e^{-0.091n})$ for $s=1.55$. 

\begin{table}[ht]
\caption{The capacity $\capa(S,[-s\i,s\i])$ for several values of $s$.}
\label{tab:strip-seg-err-i}%
\centerline{
\begin{tabular}{c|c|c|c|c|c}\hline
  $s$  & Estimated value      & Exact value        & Relative Error      & Time (sec) & Iterations \\ \hline %
$0.1$  & $1.703662933054205$ & $1.703662933054233$ & $1.7\times10^{-14}$ & 3.0  & 13\\
$0.25$ & $2.270486340549139$ & $2.270486340549105$ & $1.5\times10^{-14}$ & 2.9  & 13  \\
$0.5$  & $3.053400295538014$ & $3.053400295538072$ & $1.9\times10^{-14}$ & 2.7  & 12  \\
$1$    & $4.885188789695857$ & $4.885188789695905$ & $9.8\times10^{-15}$ & 2.5  & 10  \\
$1.5$  & $10.27204980801009$ & $10.27204980801069$ & $8.5\times10^{-14}$ & 6.9  & 32  \\
$1.55$ & $13.39253761345057$ & $13.39253761345013$ & $3.3\times10^{-14}$ & 12.6 & 56  \\
\hline
\end{tabular}}
\end{table}

\begin{example}\label{ex:ms-s}
$E=[-s,s]$ where $s>0$ is a real number.
\end{example}

In this case again, the exact value of the capacity of the condenser $(S,E)$ can be determined. 
Clearly,
$
z\mapsto\Psi^{-1}(z)=\tanh(z/2)
$
maps conformally the domain $S\setminus[-s,s]$ onto the domain $\D\setminus[-\tanh(s/2),\tanh(s/2)]$. Then, the M{\"o}bius transformation
\[
z\mapsto \frac{z+\tanh(s/2)}{1+z\,\tanh(s/2)}
\]
maps the domain $\D\setminus[-\tanh(s/2),\tanh(s/2)]$ onto the domain $\D\setminus[0,\tanh(s)]$. Hence, as in the previous example, we have
\begin{equation}\label{eq:strip-exact}
\capa(S,[-s,s]) = \capa(\D,[0,\tanh(s)]) = \frac{2\pi}{\mu(\tanh s)}.
\end{equation}

We now take $n=2^{10}$ and $r=0.2$ to compute $\capa(S,E)$ for several values of $s$. 
For $s=0.5$ and $s=2$, the relative error in the approximate values vs. $n$ are given in Figure~\ref{fig:sing-msi-err} (right). The error is $O(e^{-0.217n})$ for $s=0.5$ and $O(e^{-0.215n})$ for $s=2$.

\begin{example}\label{ex:ver}
$E=\i s+[-\i,\i]$ with $-0.55\le s\le0.55$.
\end{example}

\begin{figure}[ht] %
\centerline{
\scalebox{0.35}{\includegraphics[trim=0 0.0cm 0 0cm,clip]{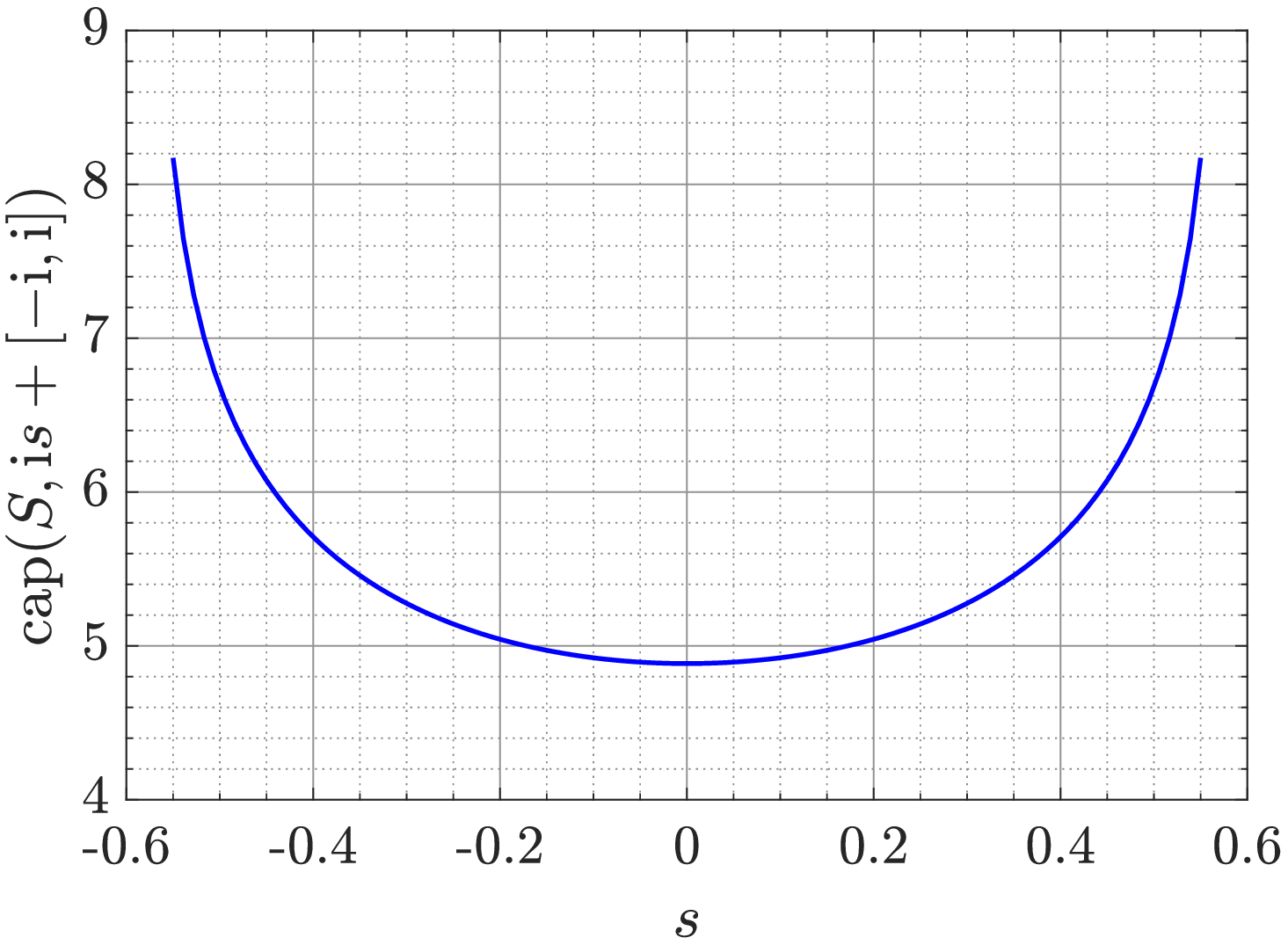}}
\hfill
\scalebox{0.35}{\includegraphics[trim=0 0.0cm 0 0cm,clip]{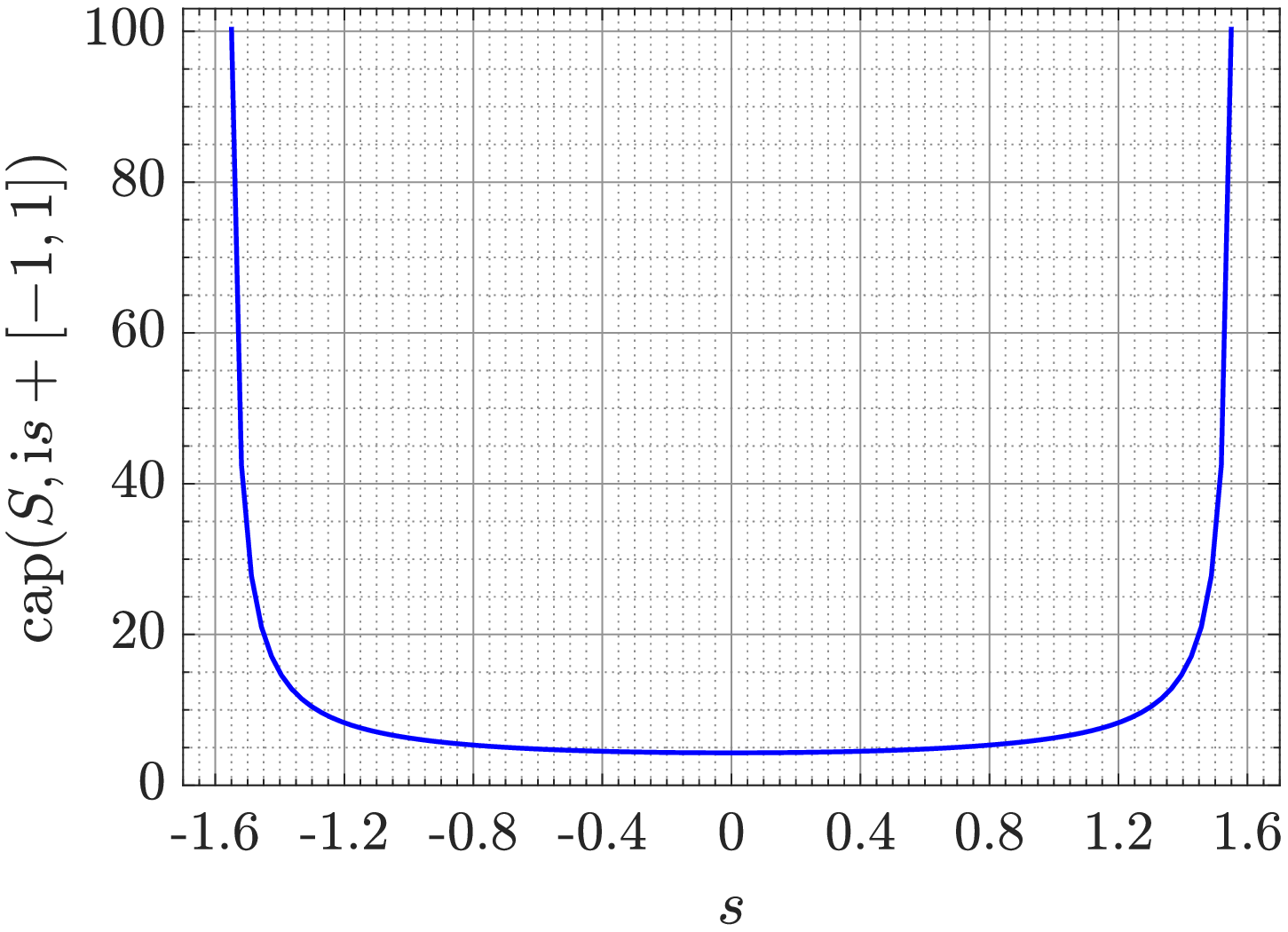}}
}
\caption{The graph of the capacity functions $\capa(S,\i s+[-\i,\i])$ in Example~\ref{ex:ver} (left) and $\capa(S,\i s+[-1,1])$ in Example~\ref{ex:hor} (right).}
\label{fig:sing-cant-mi-i}
\end{figure}

In this example, we consider the vertical segment $[-\i,\i]$ and study the effect of vertically shifting this segment on the capacity of the condenser $(S,E)$.
On the left of Figure~\ref{fig:sing-cant-mi-i}, we display the graph of the capacity $\capa(S,\i s+[-\i,\i])$ as a function of $s\in[-0.55, 0.55]$. The computation is performed with $n=2^{11}$ and $r=0.1$. 

\begin{example}\label{ex:hor}
$E=\i s+[-1,1]$ with $-1.55\le s\le 1.55$.
\end{example}

We now consider the horizontal segment $[-1,1]$ and study the effect of vertically shifting this segment on the capacity of the condenser $(S,E)$. 
The graph of the capacity as a function of $s\in[-1.55, 1.55]$, computed using the same parameters as in the previous example, is shown in Figure~\ref{fig:sing-cant-mi-i} (right).  

By looking at Figure~\ref{fig:sing-cant-mi-i}, we can notice an increase in the capacity as the segment moves vertically close the boundary of $S$. 
It is worth noting that horizontal movement of a segment in a strip does not change the capacity value because of translation invariance. In fact, for any segment $[a,b]$ in $S$ and any real number $s$, the linear transformation
\[
z\mapsto z-s
\]
maps the domain $S\backslash(s+[a,b])$ onto the domain $S\backslash[a,b]$, and hence by conformal invariance of the capacity, we have
\begin{equation}\label{eq:shit-s}
\capa(S,s+[a,b])=\capa(S,[a,b]), \quad a,b\in S, \quad s\in\R.
\end{equation}

\begin{example}
Segments with constant capacity.
\end{example}

\begin{figure}[ht] %
\centerline{
\scalebox{0.35}{\includegraphics[trim=0 0 0 0,clip]{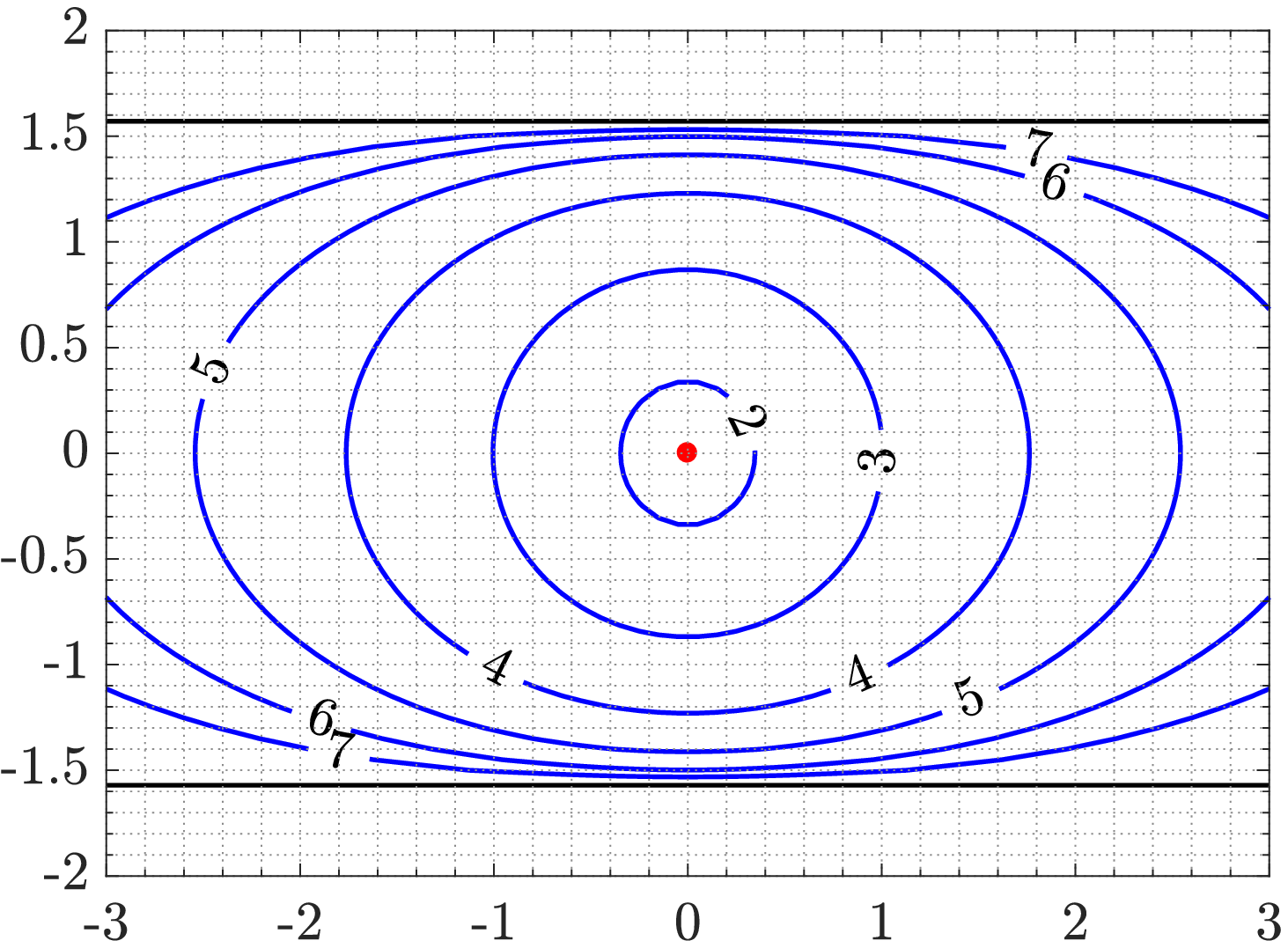}}
\hfill
\scalebox{0.35}{\includegraphics[trim=0 0 0 0,clip]{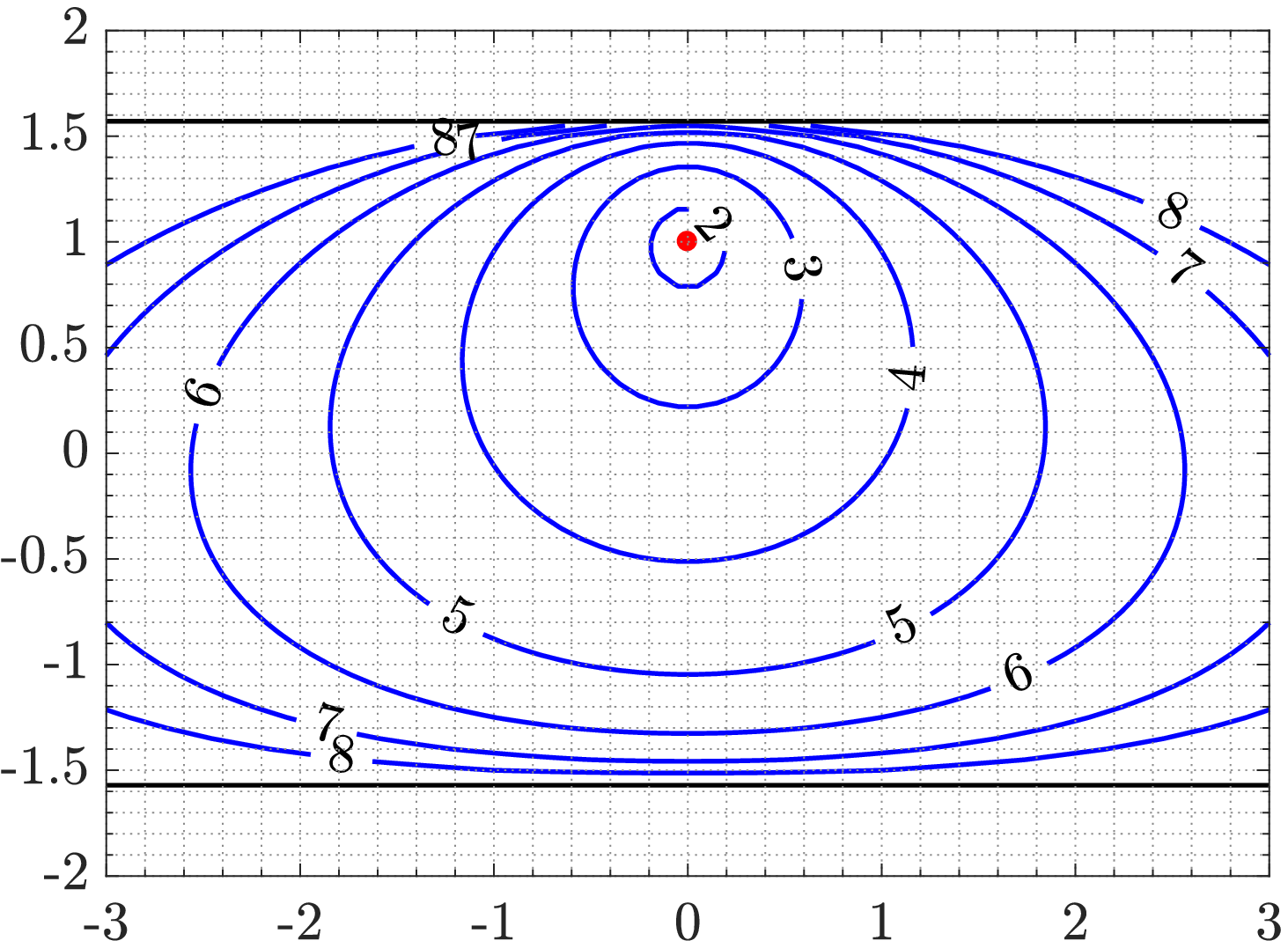}}
}
\caption{The contour lines of the capacity functions $\capa(S,[0,x+\i y])$ (left) and  the function $\capa(S,[\i,x+\i y])$ (right).}
\label{fig:sing-cant-0-yi}
\end{figure}

For a given point $a$ in $S$, consider all points $x+\i y\in S\backslash\{a\}$ such that the capacity $\capa(S,[a,x+\i y])$ is constant. 
The contour lines of the function $\capa(S,[a,x+\i y])$ in the sub-domain $x+\i y\in [-3, 3]+\i (-\pi/2, \pi/2)\setminus \{a\}$ are displayed in Figure~\ref{fig:sing-cant-0-yi} on the left for $a=0$ and on the right for $a=\i$. For each contour line, we have the same capacity for all segments with one end at $a$ and the other end on the contour line. These results are again computed with $n=2^{11}$ and $r=0.1$. 

\subsection{A strip with two slits} 

In this subsection, we want to estimate the capacity $\capa(S,E)$ where $E$ is the union of two disjoint segments.

\begin{example}
$E=[a,b]\cup[c,d]$ for several values of $a,b,c,d\, \in S$ as in Table~\ref{tab:cap-abcd}.
\end{example}

The above described method with $n=2^{10}$ and $r=0.2$ is used to compute the capacity $\capa(S,[a,b]\cup[c,d])$ for several values of $a$, $b$, $c$ and $d$. The obtained results are presented in 
Table~\ref{tab:cap-abcd}.

\begin{table}[ht]
\caption{The approximate values of the capacity $\capa(S,[a,b]\cup[c,d])$.}
\label{tab:cap-abcd}%
\centerline{
\begin{tabular}{c|c|c|c|c}\hline
$a$     & $b$     & $c$     & $d$      & $\capa(S,[a,b]\cup[c,d])$\\ \hline %
$-1$    & $-1+\i$ & $1$     & $1-\i$   & $6.0697365159628$  \\
$-1$    & $-1+\i$ & $1$     & $1+\i$   & $6.0193744425645$  \\
$-1$    & $-2$    & $1$     & $2$      & $5.6844096460738$  \\
$-1+\i$ & $1+\i$  & $-1-\i$ & $1-\i$   & $11.029565510437$  \\
\hline
\end{tabular}}
\end{table}

\begin{example}\label{ex:x+J}
$E=E_1\cup E_2$ with $E_1=-x+J$ and $E_2=x+J$ where $J=[-\i,\i]$ and $x>0$ is a real number.
\end{example}

\begin{figure}[ht] %
\centerline{\scalebox{0.4}{\includegraphics[trim=0 0 0 0,clip]{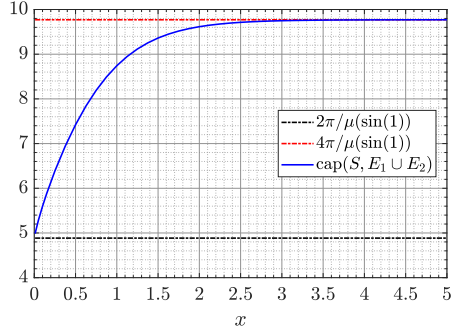}}}
\caption{Results for Example~\ref{ex:x+J}: The graph of the capacity $\capa(S,E_1\cup E_2)$ where $E_1=[-x-\i,-x+\i]$ and $E_2=[x-\i,x+\i])$.}
\label{fig:2v-u}
\end{figure}

In this example, we consider the two vertical segments $E_1$ and $E_2$ centered on the $x$-axis (middle line of the strip) and far from each other by a distance of $2x>0$. We study the effect of this distance on $\capa(S,E_1\cup E_2)$.  
We take $n=2^{11}$ and $r=\min\{0.2,x/2\}$, and compute the values of $\capa(S,[-x-\i,-x+\i]\cup[x-\i,x+\i])$ as a function of $x$ for $0.01\le x\le 4$. The graph of this function is shown in Figure~\ref{fig:2v-u}.

It follows from~\eqref{eq:shit-s} that
\[
\capa(S,-x+J)=\capa(S,x+J)=\capa(S,J).
\]
Hence, in view of~\eqref{eq:strip-exact-i}, we have
\[
\capa(S,E_1)=\capa(S,E_2)=\capa(S,[-\i,\i])=\frac{2\pi}{\mu(\sin(1))}.
\]
Then, by~\cite[Theorem 1.8]{du} and \cite[Lemma 7.1]{hkv},
\[
\frac{2\pi}{\mu(\sin(1))}=\capa(S,E_1)\le \capa(S,E_1\cup E_2) \le \capa(S,E_1)+\capa(S,E_2)=\frac{4\pi}{\mu(\sin(1))},
\]
which is validated numerically by the presented method in Figure~\ref{fig:2v-u}.
Further, Figure~\ref{fig:2v-u} suggests that 
\[
\lim_{x\to0^+}\capa(S,E_1\cup E_2)=\capa(S,[-\i,\i])=\frac{2\pi}{\mu(\sin(1))},
\]
and for large $x$,
\[
\capa(S,E_1\cup E_2)\approx \capa(S,E_1)+\capa(S,E_2)=\frac{4\pi}{\mu(\sin(1))}.
\]

\begin{example}\label{ex:x+I}
$E=E_1\cup E_2$ with $E_1=-x+I$ and $E_2=x+I$ where $I=[-1,1]$ and $x>1$ is a real number.
\end{example}

We consider here the two horizontal segments $E_1$ and $E_2$ located on the $x$-axis where the distance between them is $2x>2$. 
Note that $\lim_{x\to1^+}E_1\cup E_2=[-2,2]$. Further, it follows from~\eqref{eq:strip-exact} and~\eqref{eq:shit-s} that 
\[
\capa(S,E_1)=\capa(S,E_2)=\capa(S,I)=\frac{2\pi}{\mu(\tanh(1))}.
\]

\begin{figure}[ht] %
\centerline{\scalebox{0.4}{\includegraphics[trim=0 0 0 0,clip]{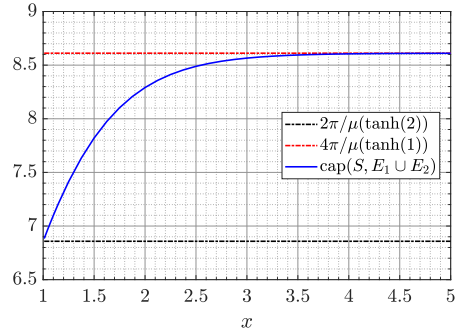}}}
\caption{Results for Example~\ref{ex:x+I}: The graph of the function $\capa(S,E_1\cup E_2)$ where $E_1=[-x-1,-x+1]$ and $E_2=[x-1,x+1]$.}
\label{fig:2h-u}
\end{figure}

We use our method with $n=2^{11}$ and $r=0.2$ to compute the values of $\capa(S,E_1\cup E_2)$ for $1.01\le x\le 4$. The graph $\capa(S,E_1\cup E_2)$, as a function of $x$, is shown in Figure~\ref{fig:2h-u}. The obtained numerical results show that
\begin{align*}
\frac{2\pi}{\mu(\tanh(2))} &=\capa(S,[-2,2])\\
						   &\le \capa(S,E_1\cup E_2) \\
						   &\le \capa(S,E_1)+\capa(S,E_2)=\frac{4\pi}{\mu(\tanh(1))}.
\end{align*}
Moreover, Figure~\ref{fig:2h-u} reveals that 
\[
\lim_{x\to1^+}\capa(S,E_1\cup E_2)=\capa(S,[-2,2])=\frac{2\pi}{\mu(\tanh(2))},
\]
and for large $x$,
\[
\capa(S,E_1\cup E_2)\approx \capa(S,E_1)+\capa(S,E_2)=\frac{4\pi}{\mu(\tanh(1))}.
\]

Figure~\ref{fig:2h-u} also leads to a new conjecture about the capacity of condensers in the strip domain $S$. Indeed, if the condenser is of the form
$(S,E)$ and $E= E_1 \cup E_2 \subset {\mathbb R}$ , $E_1$ and $E_2$ are segments, 
then
$$ {\rm cap}(S,E)\ge {\rm cap}(S,H)$$
where $H \subset {\mathbb R}$ is a segment with diameter equal to the sum of diameters
of $E_1$ and $E_2.$

\begin{example}\label{ex:xi+I}
$E=E_1\cup E_2$ where $E_1=-x\i+I$, $E_2=x\i+I$, $I=[-1,1]$, and $x>0$ is a real number.
\end{example}

\begin{figure}[ht] %
\centerline{\scalebox{0.4}{\includegraphics[trim=0 0 0 0,clip]{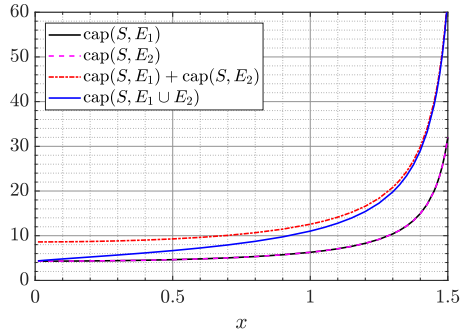}}}
\caption{Results for Example~\ref{ex:xi+I}: The graph of the function $\capa(S,E_1\cup E_2)$ where $E_1=-x\i+I$, $E_2=x\i+I$, $I=[-1,1]$, and $x>0$ is a real number.}
\label{fig:2g-u}
\end{figure}

In the final example concerning the case of two slits, we consider the two parallel horizontal segments $E_1$ and $E_2$, where both of them are centered on the $y$-axis, and study the effect of the distance between them, which is  $2x$, on the capacity $\capa(S,E_1\cup E_2)$. 
We use our method with $n=2^{12}$ and $r=\min\{0.2,x/2,1.53-x\}$ to compute the values of $\capa(S,E_1\cup E_2)$ as a function of $x$ for $0.01\le x\le 1.5$ (see Figure~\ref{fig:2g-u} for the graph of this function). For this case, the exact values of the capacities $\capa(S,E_1)$ and $\capa(S,E_1)$ are unknown and hence will be computed numerically using the proposed method.

By looking at Figure~\ref{fig:2g-u}, we can see that
\[
\capa(S,E_1)=\capa(S,E_2) \le \capa(S,E_1\cup E_2) \le \capa(S,E_1)+\capa(S,E_2).
\]
It is also clear that $\capa(S,E_1\cup E_2) \approx \capa(S,E_1)$ for small $x$, and $\capa(S,E_1\cup E_2) \approx \capa(S,E_1)+\capa(S,E_2)$ for large $x$.

\subsection{A strip with many rectilinear slits} 

The last two examples illustrate that our iterative approach can be employed in the case of a multiply connected infinite strip with a high connectivity for which the exact value of the capacity is not known.

\begin{example} We first consider the case of four rectilinear slits
$E_1=[2-\i , 3.5+0.5\i]$, $E_2=[1+\i , -1+\i]$, $E_3=[-\i , -2.5+0.5\i]$ and $E_4=[-3-\i , -3+\i]$ with $\delta_1=1$, $\delta_2=2$, $\delta_3=3$ and $\delta_4=4$.
The approximate value of $\capa(S,E,\delta)$ obtained with $n=2^{11}$ and $r=0.2$ is
\[
\capa(S,E,\delta)=41.8434999283923.
\]
\end{example}

\begin{figure}[ht] %
	\centerline{
		\scalebox{0.35}{\includegraphics[trim=0 0  0 0,clip]{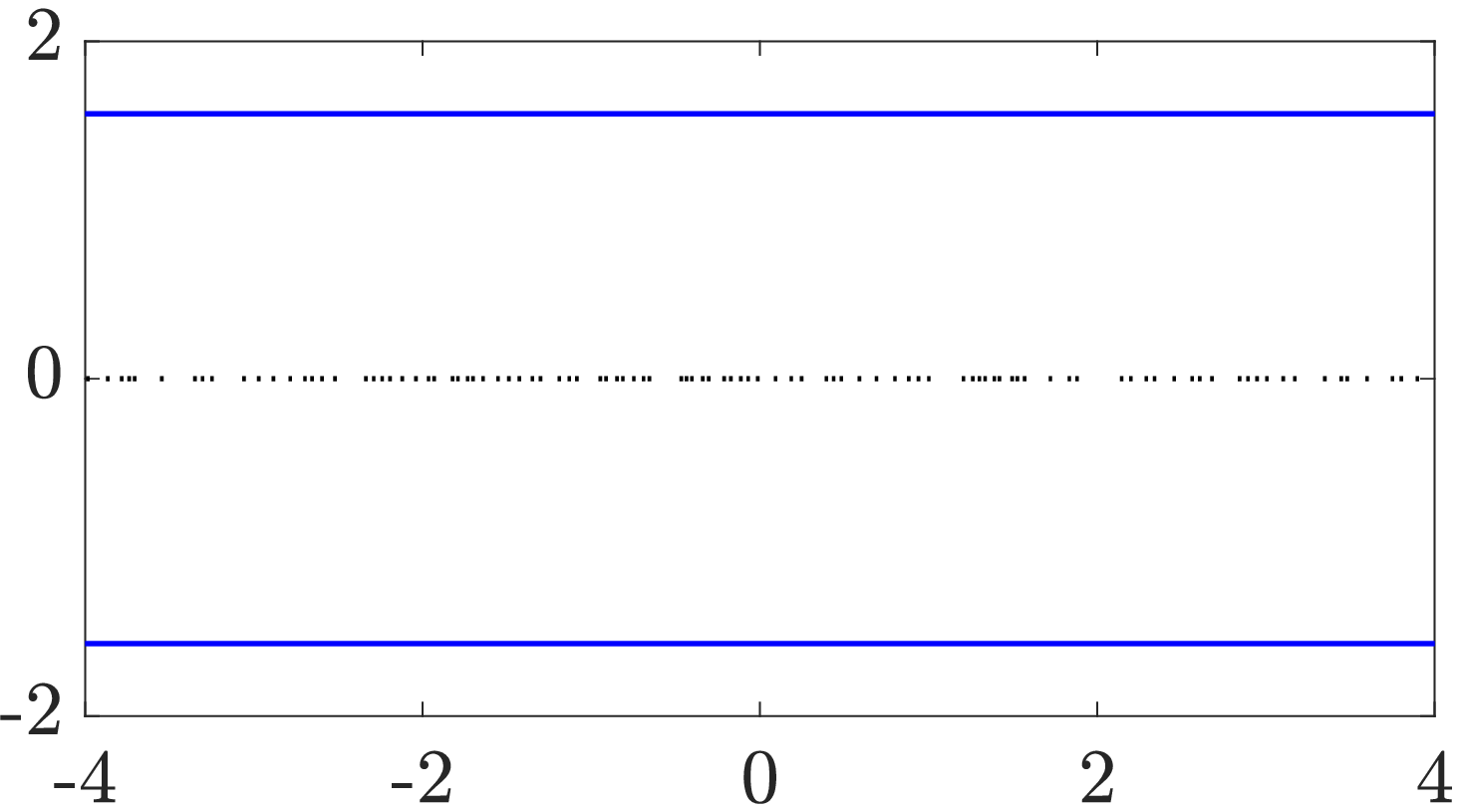}}
		\hfill
		\scalebox{0.35}{\includegraphics[trim=0 0  0 0,clip]{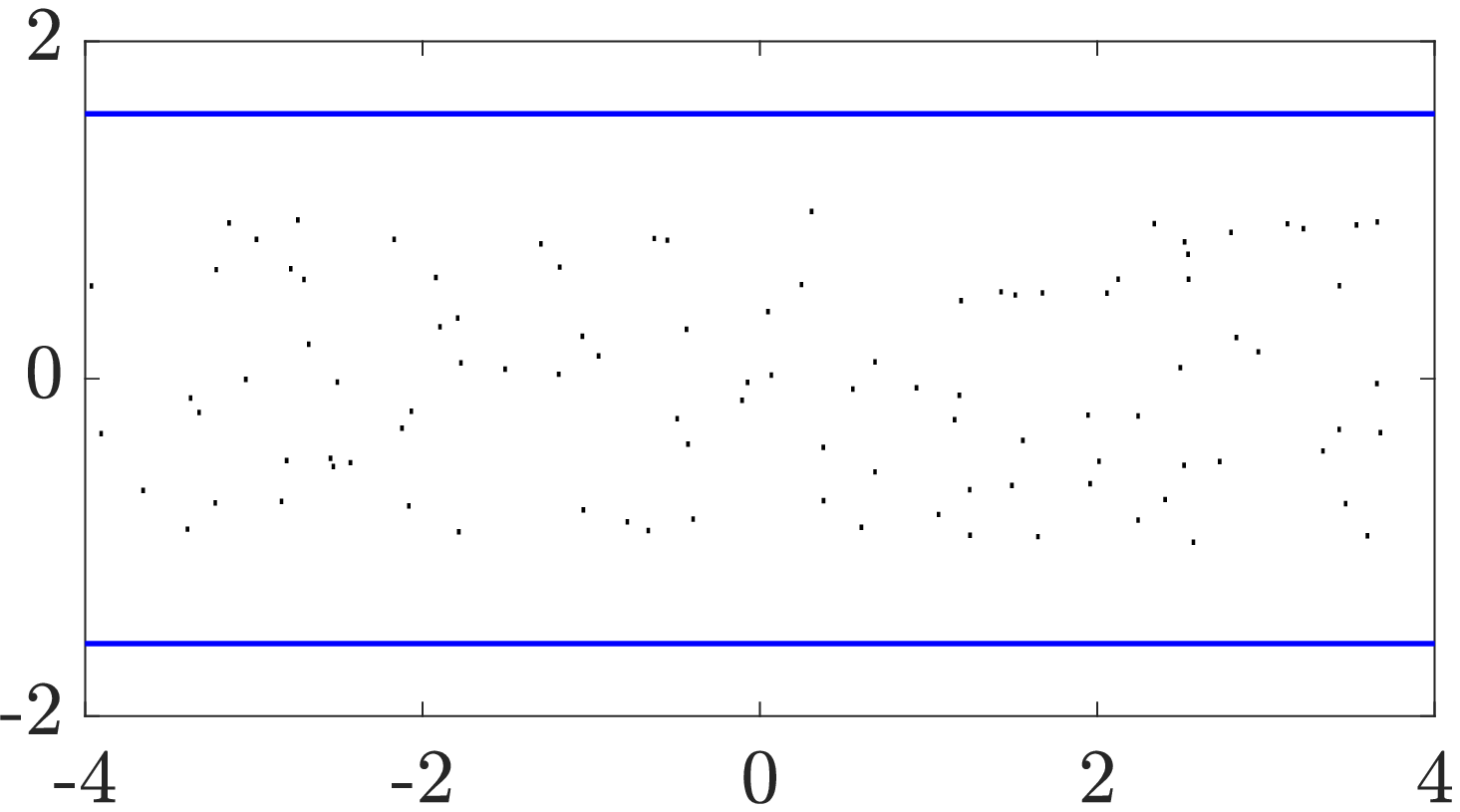}}
	}
	\caption{The domains $S\setminus \cup_{j=1}^{m} E_j$ (left) and $S\setminus \cup_{j=1}^{m}\hat E_j$ (right) in Example~\ref{ex:100} for $m=100$.}
	\label{fig:100}
\end{figure}

\begin{example}\label{ex:100}
We finally take up a collection of $m$ disjoints horizontal intervals of length $2/m$ with random location on the real axis between $-4$ and $4$; see Figure~\ref{fig:100} (left) for $m=100$.
Clearly, the sum of the diameters of all these intervals is $2$ which is equal to the diameter of $E=[-1,1]$. Note that by~\eqref{eq:strip-exact} the exact capacity involving the single plate $E$ is $\capa(S,E)=2\pi/\mu(\tanh(1))$. We  also consider the case when the $m$ collection of disjoints horizontal intervals is randomly distributed in the strip $[-4,4]\times[-1,1]$, and denote it if so by $(\hat E_j)_{j=1}^m$; see Figure~\ref{fig:100} (right) for $m=100$. 
\end{example}

\begin{figure}[ht] %
\centerline{\scalebox{0.4}{\includegraphics[trim=0 0 0 0,clip]{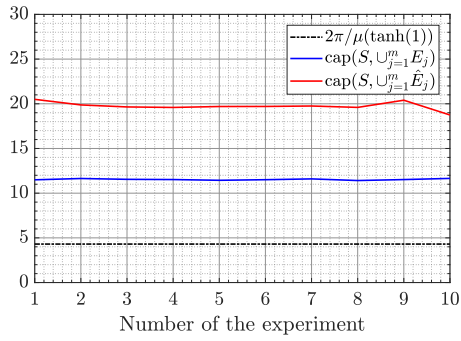}}}
\caption{Results for Example~\ref{ex:100}: The computed values of the capacities for $10$ random locations of the slits.}
\label{fig:100-10}
\end{figure}

To estimate $\capa\left(S,\cup_{j=1}^{m}E_j\right)$ and $\capa\left(S,\cup_{j=1}^{m}\hat E_j\right)$, we employ our method using $n=2^{10}$ and $r=0.2$. We run the code for $10$ times so that to get $10$ different locations for these slits. The computed values of the capacities for these $10$ locations are shown in Figure~\ref{fig:100-10}. Note that $\cup_{j=1}^{m}E_j$ and $\cup_{j=1}^{m}\hat E_j$ have the same diameter as $E$ where the center of the horizontal slits $E_j$ are real numbers in $[-4,4]$ and the centers of the horizontal slits $\hat E_j$ are complex numbers in $[-4,4]\times[-1,1]$. The curves displayed in Figure~\ref{fig:100-10} indicate that regardless of the slits positions, we have the following inequalities
\[
\capa\left(S,\cup_{j=1}^{m}\hat E_j\right) > \capa\left(S,\cup_{j=1}^{m} E_j\right)
> \capa\left(S,E\right).
\]
This result is somehow expected as the slits $\hat E_j$ are less condensed and closer to the boundary of $S$ when compared to the slits $E_j$.

\section{Uniform potential flow in multiply connected channel domains}

\begin{figure}[ht] %
	\centerline{
		\scalebox{0.3}{\includegraphics[trim=0 0 0 0,clip]{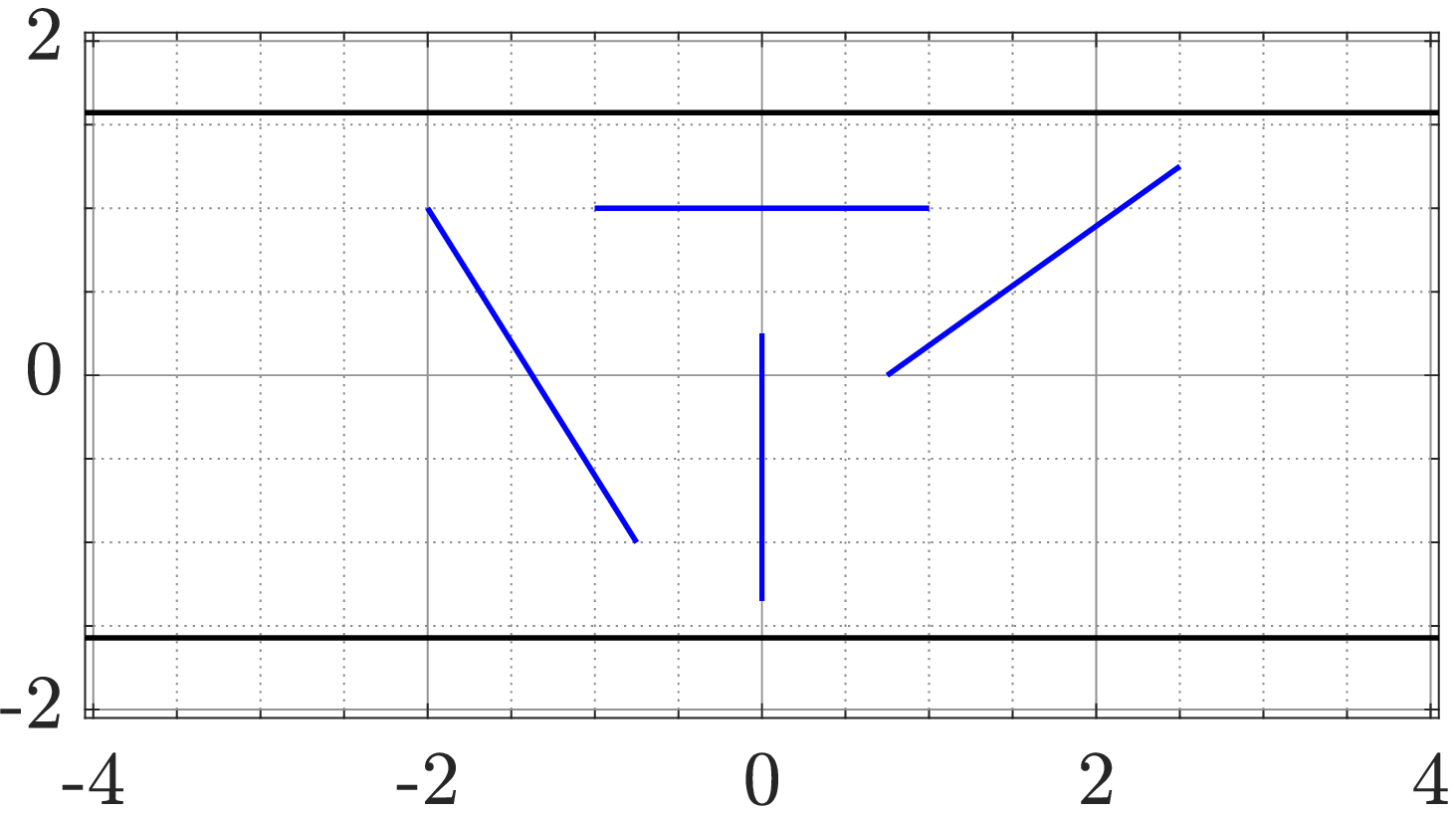}}
		\hfill
		\scalebox{0.3}{\includegraphics[trim=0 0 0 0,clip]{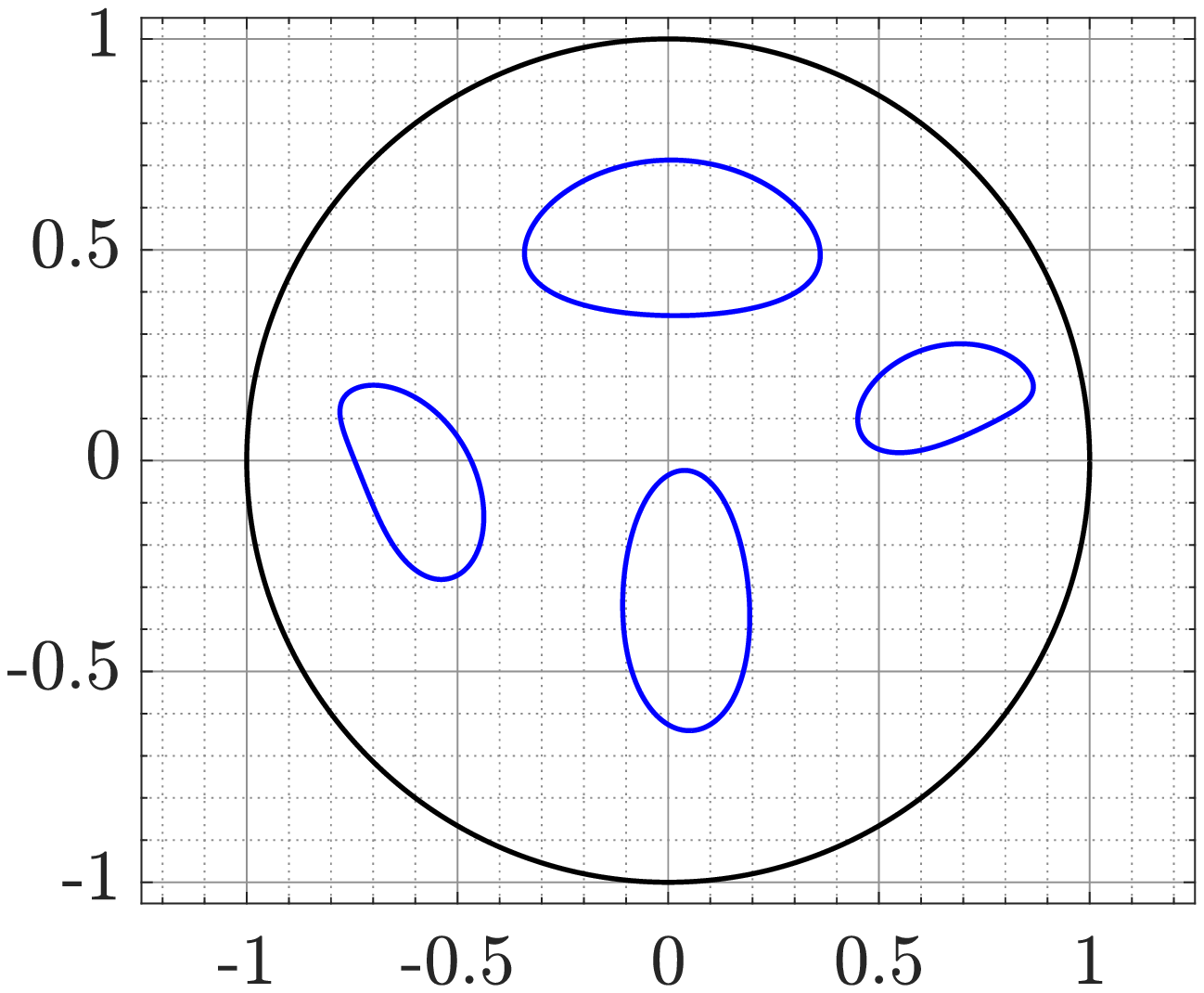}}
		\hfill
		\scalebox{0.3}{\includegraphics[trim=0 0 0 0,clip]{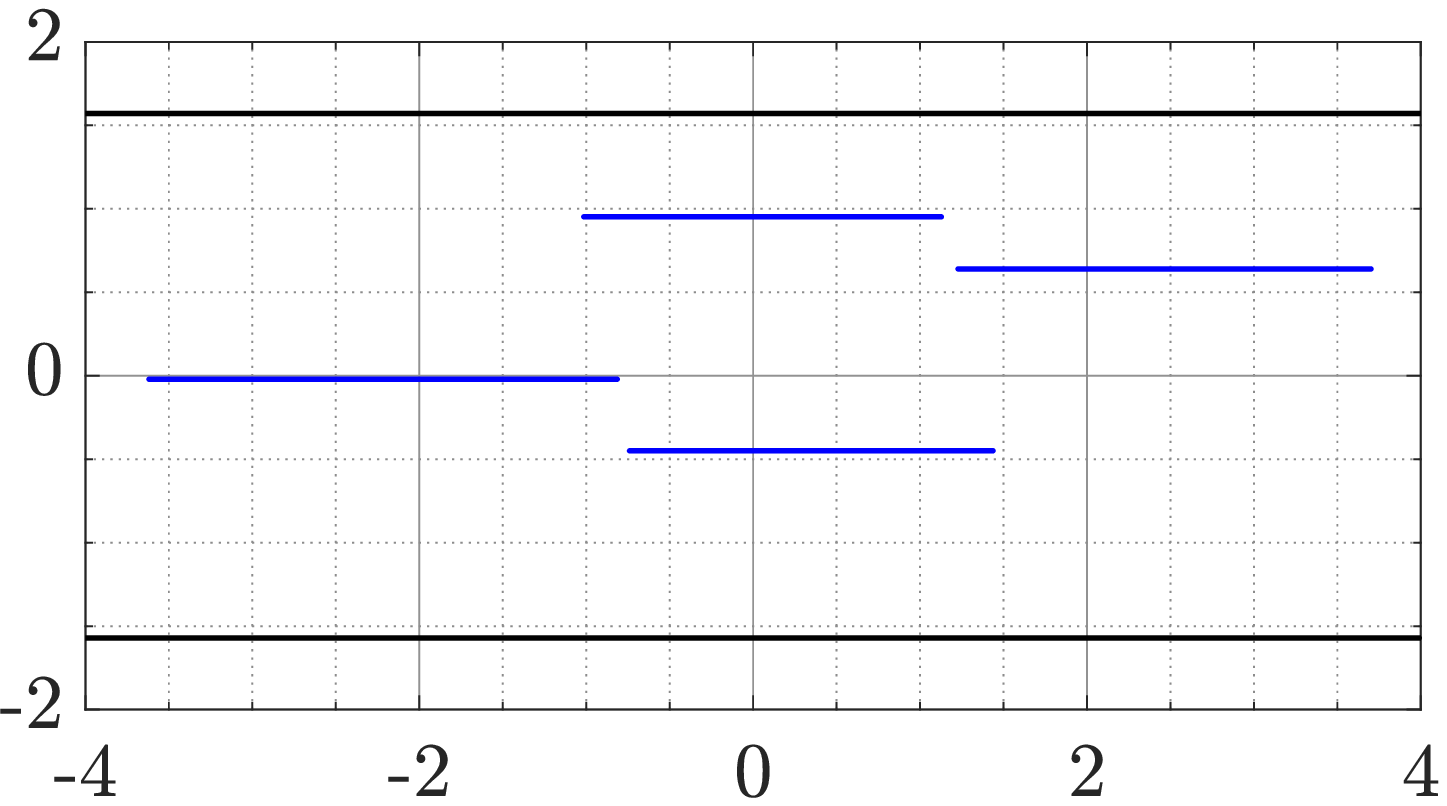}}
	}
	\caption{The domain $\Omega$ (left), the domain $G$ (center), and the domain $H$ (right) for $m=4$. The domain $G$ is computed with $n=2^{11}$ and $r=0.5$.}
	\label{fig:channel-flow}
\end{figure}

In this section, we present a fast and accurate numerical method for constructing the complex potential function $W(z)$ for a uniform incompressible, inviscid and irrotational flow past multiple disjoint segment obstacles in the strip $S=\{z\,:\, |\Im z|<-\pi/2\}$ in the case when the circulations around the segments are zeros. 
The method is based on using the above iterative method to construct a conformal mapping $w=F(z)$ from a domain $\Omega$ obtained by removing $m$ non-overlapping rectilinear slits from the strip $S$ (see Figure~\ref{fig:channel-flow} (left)) onto a domain $H$ obtained by removing $m$ horizontal slits from the strip $S$ (see Figure~\ref{fig:channel-flow}(right)). Then 
\[
W(z)=F(z)
\] 
represents the complex potential for a uniform flow in the domain $\Omega$ and the level curves of $\Im[W(z)]$ represent streamlines of the flow.

\begin{figure}[ht] %
\centerline{
\scalebox{0.4}{\includegraphics[trim=0 2cm 0 2cm,clip]{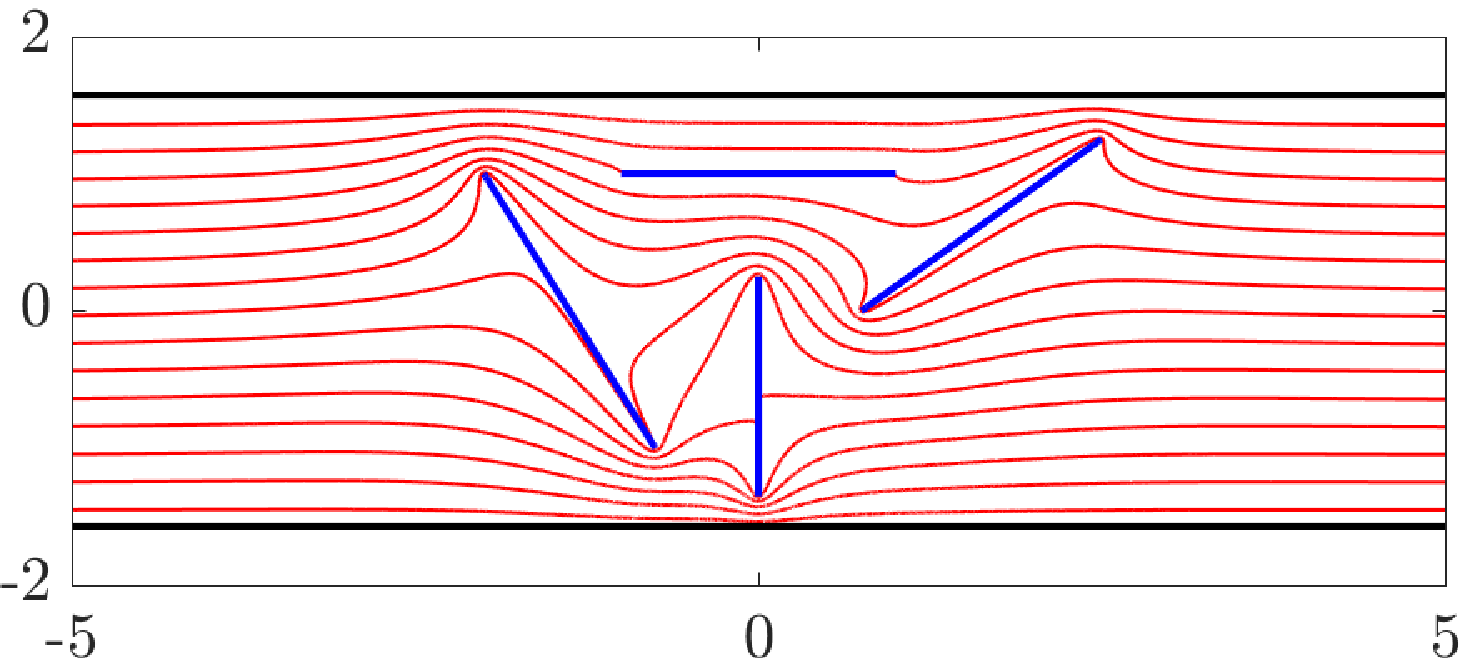}}
\hfill
\scalebox{0.4}{\includegraphics[trim=0 2cm 0 2cm,clip]{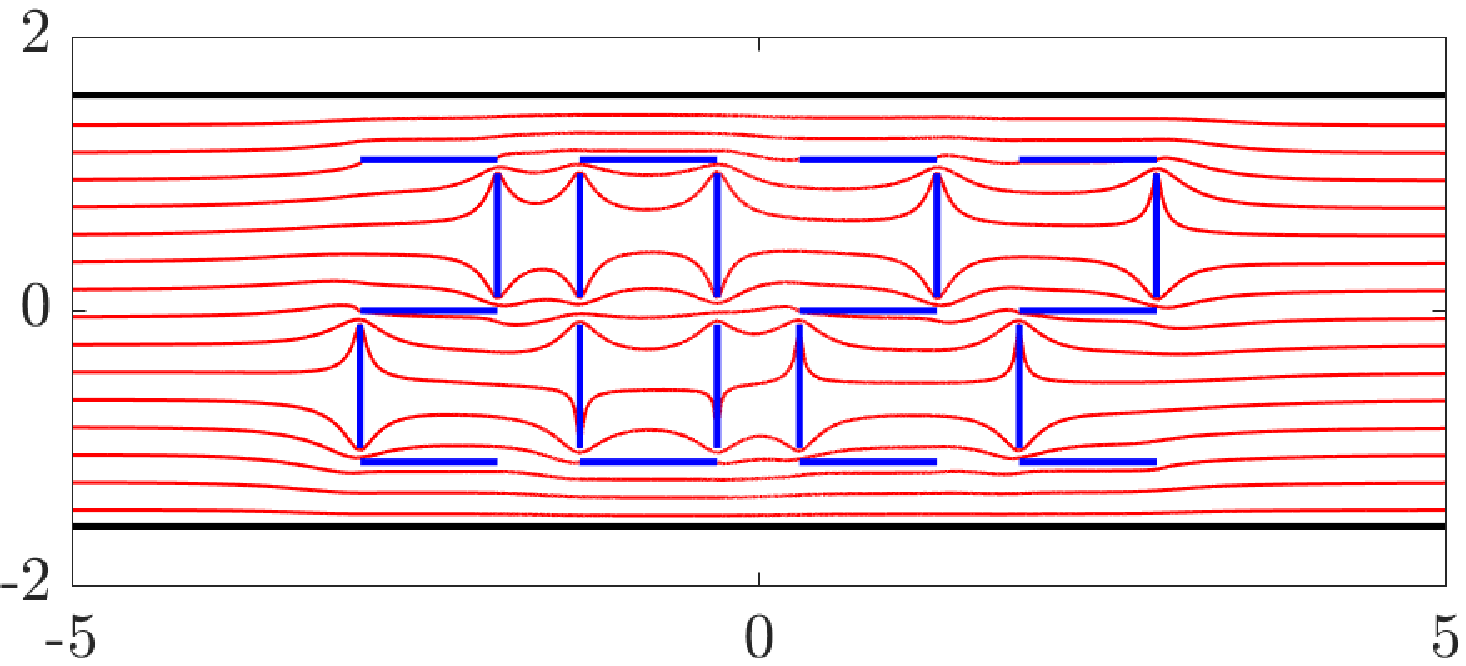}}
}
\caption{Streamlines for a uniform flow in the multiply connected channel domain $\Omega$ for $m=4$ (left) and $m=21$ (right).}
\label{fig:channel-uflow}
\end{figure}

To construct the conformal mapping $w=F(z)$, we first use the iterative method described in Section~\ref{sec:strip} to obtain a preimage domain $G$ bordered by smooth Jordan curves and the conformal mapping $z=\Psi(\zeta)$ from $G$ onto $\Omega$ (see Figure~\ref{fig:channel-flow}(center)). Then the method presented in Theorem~\ref{thm:cm-strip} is used to compute a conformal mapping $w=\Upsilon(\zeta)$ from the domain $G$ on the domain $H$. Thus, the function
\[
F(z)=(\Upsilon\circ\Psi^{-1})(z)
\]
is the required conformal mapping from the given domain $\Omega$ onto the domain $H$.
We have applied the method with $n=2^{11}$ and $r=0.2$ to compute the streamlines of a uniform flow for two channel domains as shown in Figure~\ref{fig:channel-uflow}.

\section{Conclusion} 	

In this paper, for a given the multiply connected slit domain 
\[
\Omega=S \setminus \cup_{j=1}^{m}[a_j,b_j]
\]
where $S=\{z\,:\, \left|\Im z\right|<\pi/2\}$ and $a_j,b_j\in S$ are complex numbers, $j=1,\ldots,m$, we presented an iterative method for computing a preimage multiply connected domain $G$ bordered by smooth Jordan curves as well as a conformal mapping from $G$ onto $\Omega$. 

The proposed iterative method will be useful in applications that require solving boundary value problems in the domain $\Omega$. Assuming that the considered boundary value problem is invariant under conformal mappings, then solving the problem will be much easier in the domain $G$ with smooth boundaries compared to the original domain $\Omega$ bordered by slits. 
Two applications were presented in this paper to illustrate the performance of the proposed method. 
The first application was about computing the capacity of condensers $(S,E)$ where $E\subset S$ is a union of disjoint segments. 
In the second application, we computed the streamlines for uniform flow past disjoint segments in the strip $S$. 
The presented examples demonstrated the effectiveness of the proposed iterative method even when the slits are close to each other and for domains of high connectivity. 
	
Finally, to ensure the convergence of the iterative method when the slits are close to each other, we need to choose small values of the parameter $r$ (see also~\cite{NG18,N19}). However, for given domains $\Omega$ with well-separated slits, we can choose $r=1$ and hence the inner curves in the intermediate domain $\hat\Omega^0$ and consequently the inner curves in the preimage domain $G$ will be circles, which means that $G$ is a circular domain (see Figure~\ref{fig:strip}). 
Note that, for circular multiply connected domains, analytic formulas for several fluid problems have been recently derived by Crowdy in terms of a special transcendental function known as the Schottky-Klein prime function, see~\cite{crow-10,crow-20} and the references cited therein. 
A key feature of our proposed iterative method is that in mapping the slit domain $\Omega$ onto a circular domain $G$, we can use these analytic formulas to solve several fluid problems in the slit domain $\Omega$.

\begin{acknowledgements}
We are indebted to Prof. A. Yu. Solynin  and Prof. D. Betsakos who have independently provided an analytic argument to confirm our experimental discovery~\eqref{eq:ineq} and Example~\ref{ex:x+J}. We would also like to thank two anonymous reviewers for their valuable comments and for bringing several bibliographic items to our attention. 
\end{acknowledgements}

%
\section*{Conflict of interest}
 
The authors declare that they have no conflict of interest.




\end{document}